\documentclass[twoside,11pt]{article}

\usepackage{amsmath,amssymb,amsthm}
\usepackage{mathtools}

\usepackage{blindtext}
\usepackage{color}

\usepackage{jmlr2e}

\theoremstyle{definition}
\newtheorem{theo}{Theorem}[section]

\newtheorem{prop}[theo]{Proposition}

\newtheorem{defn}[theo]{Definition}

\newtheorem{rem}[theo]{Remark}
\newtheorem{assumption}[theo]{Assumption}

\usepackage{lastpage}
\jmlrheading{}{}{}{}{}{}{}

\ShortHeadings{}{}
\firstpageno{1}

\begin{document}

\title{Exponential Inequalities for Some Mixing Processes and Dynamic Systems}

\author{\name Zihao Yuan \email Zihao.Yuan@ruhr-uni-bochum.de \\
       \addr Faculty of Mathematics\\
       Ruhr University Bochum\\
       Universitätsstraße 150, D-44780, Bochum
       \AND
       \name Holger Dette \email Holger.Dette@ruhr-uni-bochum.de \\
       \addr Faculty of Mathematics\\
       Ruhr University Bochum\\
       Universitätsstraße 150, D-44780, Bochum}

\editor{}

\maketitle

\begin{abstract}

Many important dynamic systems, time series models or even algorithms exhibit non-strong mixing properties. In this paper, we introduce the general concept of  $\mathcal{C}_{p,\mathcal{F}}$-mixing to cover such cases, where assumptions on the dependence structure become stronger with increasing $p\in [1, \infty].$
We derive  a series of sharp exponential-type (or Bernstein-type) inequalities under this dependence concept for $p=1$ and $p=\infty$. 
 More specifically,  $\mathcal{C}_{\infty,\mathcal{F}}$-mixing is equal to the widely discussed $\mathcal{C}$-mixing \citep{maume2006exponential}, and we prove a refinement  of an Berntsein-type inequality in  \cite{hang2017bernstein} for  $\mathcal{C}$-mixing processes  under more general assumptions. As there exist many stochastic processes and dynamic systems, which are not  $\mathcal{C}$ (or $\mathcal{C}_{\infty,\mathcal{F}}$)-mixing, we derive  Bernstein-type inequalities for 
 $\mathcal{C}_{1,\mathcal{F}}$-mixing processes as well and we use this result to investigate the convergence rates of plug-in-type estimators of the 
 local conditional mode set for vector-valued output, in particular in situations where the density is less smooth. 
 
\end{abstract}


\parindent 0cm 

\section{Introduction} \label{sec1} 
\def\theequation{1.\arabic{equation}}	
   \setcounter{equation}{0}

Concentration inequalities provide finite sample bounds on the deviation of an estimator from its expectation, which is crucial for understanding of the stability and reliability of statistical methods and machine learning algorithms. Hoeffding’s, Bernstein’s, McDiarmid’s or Talagrand’s inequalities are key techniques in proving non-asymptotic  statistical guarantees, such as convergence rates of estimators and hypothesis testing procedures, and we refer to the monographs  of  \cite{devroye2001combinatorial}, \cite{gyorfi2002distribution}, \cite{devroye2013probabilistic} for more details and the references therein. However, a large part of the literature  considers such inequalities under the assumption of independent data,
which is often not satisfied in areas such as spatial data regression, time series forecasting, and text or speech recognition. Consequently, numerous researchers have worked on the development of concentration inequalities for dependent data. 
Among others, we mention the work of \cite{bosq1993bernstein,modhamasry1996,chen1998sieve,merlevede2011bernstein} who derived 
and applied  Bernstein-type inequalities for 
empirical-mean-type statistics$\frac{1}{n}\sum_{i=1}^{n}h(Z_i)$, where $(Z_i)_{i\in \mathbb{N}}$ is a real valued  $\alpha$- (strong-), $\beta$- mixing or $\phi$-mixing process and $h$ denotes a real-valued measurable function.
\\
  Strong- ($\alpha$-)mixing is the weakest among these three mixing concepts, but as pointed out by many authors (see \cite{rio2017asymptotic}, \cite{doukhan2012mixing}, among others),
  there exist many important time series models and dynamic systems, which are not strong-mixing, and several authors have proposed weaker mixing concepts to address these processes. For example, \cite{rio1996theoreme}  and \cite{olivier2010deviation}  introduced generalizations of $\phi$-mixing  and \cite{dedecker2007weak} developed mixing coefficients extending $\beta$- and $\alpha$- mixing. As these concepts of dependence modeling are still not including many of the commonly considered dynamical systems ,  \cite{maume2006exponential} introduced 
   $\mathcal{C}$-mixing.
   Recently,  \cite{hang2017bernstein}  
   established  a Bernstein-type inequality for geometrically $\mathcal{C}$-mixing processes and  used this result to prove that several  state-of-the-art learning methods  can recover the (essentially) optimal rates  from  the iid case.

 In this paper we contribute to this literature from several perspectives.  First, 
in Section \ref{sec2}, we introduce the very general  concept of $\mathcal{C}_{p,\mathcal{F}}$-mixing (here $p\in [1,\infty]$ is  parameter and ${\cal F}$ is a class of functions defining different mixing coefficients), which incorporates most of the commonly used  mixing conditions for  stochastic processes and dynamic systems that are not strong mixing. In particular, our most restrictive case, $\mathcal{C}_{\infty,\mathcal{F}}$-mixing is equal to the widely used $\mathcal{C}$-mixing condition \citep{maume2006exponential}. Second, under this dependence structure, 
we  derive  sharp exponential inequalities for the empirical mean $\frac{1}{n}\sum_{i=1}^{n}h(Z_i)$, where $h$ is a real-valued measurable function (which can depend on $n$) and $(Z_i)_{i \in \mathbb{N}}$ is a stochastic process taking value in measurable space $(\mathbf{Z},\mathcal{Z})$. More specifically, 
for geometrically  $\mathcal{C}_{\infty,\mathcal{F}}$-mixing processes we derive a sharper 
 Bernstein-type inequality. Compared to Theorem 3.1 in \cite{hang2017bernstein}, our inequality is applicable  under even more general assumptions (see Section \ref{sec31}). Third, we also take a step forward by showing a sharp exponential inequality for the empirical mean whose data is collected from  a $\mathcal{C}_{1,\mathcal{F}}$-mixing process 
 taking values in a general measurable space  (see Section \ref{sec32}). Deriving such results is challenging as,  on one hand, common techniques to bound the Laplace  transform  are 
 not applicable, and, on the other hand, many of the available coupling lemmas can not be used, since we are interested in an inequality for processes taking values  in a general measurable space.
Fourth, in Section \ref{sec4}, we demonstrate the applicability of our results and  derive  the convergence rate of an estimator of  the  local conditional mode set of a conditional density with respect to Hausdorff loss. Unlike \cite{chen2016nonparametric}, which concentrate  on  iid processes with a very smooth joint density, we also consider  the scenario where the  joint density is less smooth  at its conditional mode points and where the data comes from a  geometrically   $\mathcal{C}_{1,\mathcal{F}}$-mixing process.
 
\section{$\mathcal{C}_{p,\mathcal{F}}$-mixing processes}
\label{sec2}
\def\theequation{2.\arabic{equation}}	
   \setcounter{equation}{0}
In this section, we introduce the concept of “$C_{p,\mathcal{F}}$-mixing”  and illustrate its significant generality. Throughout this paper, 
 we consider random variables defined on  a common probability space $(\Omega, \mathcal{A},\mathbb{P})$   taking  values in  a measurable space $(\mathbf{Z},\mathcal{Z})$ and the  class  
\begin{align} 
\label{det1}
\mathcal{F}=\{f: \mathbf{Z}: \to \mathbb{R} ~|~ f \text{ is a measurable such that}\  ||f||< \infty\},
\end{align} 
where $\| \cdot \|$ is a given (semi-)norm. Note that the (semi-)norm defines the class ${\cal F} $ although we do not reflect this in our notation. Similarly, we will consider later functions $ f=f_n \in {\cal F} $ which depend on the sample size without making this explicit in the notation. 

\begin{defn}
    \label{def 1} 
    \it {Let $p, q\in [1,\infty]$ such that ${1\over p} + {1\over q} =1$.
    A $\mathbf{Z}$-valued stochastic process  
    $(Z_{n})_{n\in\mathbb{N}}$  
    is called  $\mathcal{C}_{p,\mathcal{F}}$-mixing if there exists a decreasing
    function $C: \mathbb{R} \to [0,\infty) $  
    such that the inequality
    \begin{align}
        \label{eq 2-1}
        \max_{n\in\mathbb{N}}|Cov(Y,f(Z_{n+k}))|\leq C(k)||Y||_{q}||f||
    \end{align}
    holds for any $f\in \mathcal{F}$ and   any 
    real-valued random variable $Y$ which is  measurable with respect to the sigma field $\mathcal{M}_{1}^{n}= \sigma (Z_1, \ldots , Z_n)$.  $C(k)$ is called  the ($k$-th) {$\mathcal{C}_{p,\mathcal{F}}$-mixing coefficient}. 
   $(Z_{n})_{n\in\mathbb{N}}$ is called {geometrically  } 
   $\mathcal{C}_{p,\mathcal{F}}$-mixing  of order $(\nu,b,\gamma)$, if $C(k)\leq \nu^{-bk^{\gamma}}$ for some $\nu>1$, $b,\gamma>0$, and   {algebraically}  $\mathcal{C}_{p,\mathcal{F}}$-mixing  of order $\gamma$, if 
   its mixing coefficients satisfy $C(k)\leq k^{-\gamma}$ for some $\gamma>0$.}
\end{defn}
It is quite clear that a $C_{p,\mathcal{F}}$-mixing process is also  $C_{p',\mathcal{F}}$-mixing for all $p'\geq p$. 
Therefore, the most restrictive condition is $\mathcal{C}_{\infty,\mathcal{F}}$-mixing which coincides with the definition of $\mathcal{C}$-mixing process introduced in \cite{maume2006exponential} and \cite{hang2017bernstein}.
This mixing concept has turned out to be  general enough to include many important non-strong mixing dynamic systems as special cases and we  refer to \cite{dedecker2005new} and \cite{hang2017bernstein} for some review of specific examples.  
In order to connect  Definition \ref{def 1}  with some widely used mixing concepts, we first derive a sufficient condition for a process to be  $\mathcal{C}_{p,\mathcal{F}}$-mixing. For this purpose we consider  for $B>0$ the “ball” 
 $\mathcal{F}_{B}=\{f\in\mathcal{F}:||f||\leq B\} $ and introduce the coefficients
\begin{align}
    \label{eq 2-3}
C_{p}^{*}(k):= \max_{n}\sup_{f\in\mathcal{F}_{B}} \big \| E[f(Z_{n+k})|\mathcal{M}_{1}^{n}]-E[f(Z_{n+k})] \big \|_{p}.
    \end{align}
Then, we immediately obtain the following result. 
\begin{prop}
    \label{prop 2}
    \it {If the coefficient $C_p^*(k)$ of the $\mathbf{Z}$-valued stochastic process $(Z_{n})_{n\in\mathbb{N}}$ converges to $0$, $(Z_{n})_{n\in\mathbb{N}}$ is $\mathcal{C}_{p,\mathcal{F}}$-mixing, and for any $f\in\mathcal{F}$, the inequality
    \begin{align*}
       \max_{n}|\text{Cov}(Y,f(Z_{n+k}))|\leq C(k)||Y||_{q}||f||
    \end{align*}
   holds, where ${1 \over p} + {1 \over q} =1$ and $ C(k)=B^{-1}C_{p}^{*}(k)$.}
\end{prop}
Proposition \ref{prop 2} gives a sufficient condition for a process to be $\mathcal{C}_{p,\mathcal{F}}$-mixing.  Note that the coefficients in \eqref{eq 2-3} are essentially 
built on the $L_{p}$-norm of the conditional expectation of  specific centered random variables, with respect to some sigma algebras.
More precisely, with the notation 
\begin{align}
    \label{eq 2-4}
C_{p,\mathcal{F}_{B}}(\mathcal{M},Z):=     \sup_{f\in\mathcal{F}_{B}} \big \| E[f(Z)|\mathcal{M}]-E[f(Z)] \big \|_{p},
\end{align}
we have 
\begin{align}
\label{det11}
C_{p}^{*}(k):= \max_{n} C_{p,\mathcal{F}_{B}}({\mathcal{M}_1^n,Z_{n+k}}).
\end{align}
With this point of view,  we can investigate the relation between $\mathcal{C}_{p,\mathcal{F}}$- and other 
mixing concepts. In fact we  will demonstrate that many mixing coefficients proposed in the literature are 
special cases of the coefficient \eqref{eq 2-4}.

For this purpose we first review four important mixing coefficients which have been introduced for the case $\mathbf{Z}=\mathbb{R}$ by \cite{dedecker2005new} and are defined as follows
\begin{align}
\label{equation 1}
    &\tau(\mathcal{M},Z):=\int \big \|F_{Z|\mathcal{M}}(z)-F_{Z}(z)\big \|_{1}dz,\\
\label{equation 2}
    &\Tilde{\alpha}(\mathcal{M},Z):=\sup_{z\in\mathbb{R}}\big \|F_{Z|\mathcal{M}}(z)-F_{Z}(z)\big \|_{1},\\
\label{equation 3}
    &\Tilde{\beta}(\mathcal{M},Z):=\big \|\sup_{z\in\mathbb{R}}|F_{Z|\mathcal{M}}(z)-F_{Z}(z)\big \|_{1},\\
\label{equation 4}
    &\Tilde{\phi}(\mathcal{M},Z):=\sup_{z\in\mathbb{R}}\big \|F_{Z|\mathcal{M}}(z)-F_{Z}(z)\big \|_{\infty},
\end{align}
where $F_{Z|\mathcal{M}}(z) =\mathbb{P}(Z\leq z|\mathcal{M})$ is the conditional distribution of $Z$ given $\mathcal{M}$. With these notations, we can define mixing coefficients $\tau(\mathcal{M}_{1}^n,Z_{n+k})$, $\tilde{\alpha}(\mathcal{M}_{1}^n,Z_{n+k})$, $\tilde{\beta}(\mathcal{M}_{1}^n,Z_{n+k})$ and $\tilde{\phi}(\mathcal{M}_{1}^n,Z_{n+k})$, which can be compared with \eqref{det11}.
Note that we use the notation $\tilde \alpha$, $\tilde \beta$ and $\tilde \phi$ to distinguish these coefficients from those appearing in the definition of  $\alpha$-(strong-), $\beta$- and $\phi$-mixing \citep[see, for example,][]{doukhan2012mixing} and
that  $\Tilde{\alpha}(\mathcal{M},Z)\leq \Tilde{\beta}(\mathcal{M},Z)\leq \Tilde{\phi}(\mathcal{M},Z)$. 
As pointed out in 
\cite{dedecker2005new} and \cite{dedecker2007weak},
the coefficients \eqref{equation 1}-\eqref{equation 4} cover many non-strong mixing time series or dynamic systems. To explain how they are included in Definition \ref{def 1}, we note that, by  Lemma 1 in \cite{dedecker2005new}, 
there exist alternative representations for these coefficients that are
\begin{align}
\label{eq 2-5}
    & \tau (\mathcal{M},Z)=\Big \| \sup\Big \{\Big |\int f d\mathbb{P}_{Z|\mathcal{M}}(x)-\int f d\mathbb{P}_{Z}(x)\Big |, f\in \text{Lip}_{1}(\mathbb{R})\Big \}\Big \|_{1},\\
\label{alpha}
    & \Tilde{\alpha}(\mathcal{M},Z)=\sup\big \{ \|E[f(Z)|\mathcal{M}]-E[f(Z)] \|_{1}:f\in \text{BV}_{1}(\mathbb{R}) \big \},\\
\label{beta}
    & \Tilde{\beta}(\mathcal{M},Z)=\Big \| \sup\Big \{\Big |\int f d\mathbb{P}_{Z|\mathcal{M}}(x)-\int f d\mathbb{P}_{Z}(x)\Big | f\in \text{BV}_{1}(\mathbb{R})\Big \}\Big \|_{1},\\
\label{eq 2-6}
    & \Tilde{\phi}(\mathcal{M},Z)=\sup \big \{ \|E[f(Z)|\mathcal{M}]-E[f(Z)] \|_{\infty} \big |~f\in \text{BV}_{1}(\mathbb{R}) \big \}.
\end{align}
Here $\mathbb{P}_{Z|\mathcal{M}}$ is the conditional  distribution 
of the random variable $Z$  with respect to $\sigma $-field ${\cal M}$ and 
$\text{Lip}_{1}(\mathbb{R})$ and $\text{BV}_{1}(\mathbb{R})$ are  the unit balls in the  normed spaces of functions $f : \mathbb{R} \to \mathbb{R}$ equipped with the  Lipschitz-norm  $|f|_{\text{Lip}}$ and  the  total variation (\textbf{TV})-norm  $|f|_{TV}$, respectively.  Apparently, the coefficients $C_{p,\mathcal{F}_{B}}(\mathcal{M},Z)$ defined in \eqref{eq 2-4} can be regarded as a direct generalization of these four mixing coefficients. 
Due to Proposition \ref{prop 2}, it instantly becomes a special case of our $\mathcal{C}_{1,\mathcal{F}}$-mixing processes defined in Definition \ref{def 1}.

We emphasize that $\mathcal{C}_{1,\mathcal{F}}$-mixing also covers the traditional $\alpha$-mixing for real valued processes. To see this we note, 
for real-valued random variables, one can show that
the representation \eqref{alpha} is equivalent to  
\begin{align*}
\tilde\alpha(\mathcal{M},Z) & =\sup\{\mathbb{P}(Z \leq z |\mathcal{M})-\mathbb{P}(Z \leq z ) ~|~z\in\mathbb{R} \} \\
& \leq  
    \alpha(\mathcal{M},\sigma(Z)) :=\sup\{\mathbb{P}(Z\in A|\mathcal{M})-\mathbb{P}(Z\in A):A\in\sigma(Z)\}  .
\end{align*}
 Hence, a $\mathcal{C}_{1,\mathcal{F}}$-mixing condition incorporates an $\alpha$-mixing condition  defined on a sub-sigma algebra and a sigma-algebra generated by one real-valued random variable. 
 Consequently, a real-valued   $\alpha$-mixing  process \citep[see][]{doukhan2012mixing}
 is  also  $\mathcal{C}_{1,\mathcal{F}}$-mixing. Moreover, as pointed out in
 Definition 2.5 and Lemma 2.1 in \cite{dedecker2007weak}, these observations are also correct in the case $\mathbf{Z}=\mathbb{R}^{d}$. Furthermore,   
 $\tau$-mixing processes taking values in a Polish space are $\mathcal{C}_{1,\mathcal{F}}$-mixing as well. This follows from  Proposition \ref{prop 2}, the inequality 
 \begin{align}
\label{equation 6}
    \tau(\mathcal{M},Z)&\geq \sup\left\{\left\|\int f d\mathbb{P}_{Z|\mathcal{M}}(x)-\int f d\mathbb{P}_{Z}(x)\right\|_{1}, f\in \text{Lip}_{1}(\mathbb{R})\right\}\notag\\
    &=\sup\{||E[f(Z)|\mathcal{M}]-E[f(Z)]||_{1}:f\in \text{Lip}_{1}\}
\end{align}
and the fact that  $\tau(\mathcal{M},Z)$ can be easily be defined for processes taking values a Polish space \citep[see][]{dedecker2005new}.  
 Summarizing this discussion,  besides the classical mixing concepts,  Definition \ref{def 1}   covers many useful non-strong mixing processes such as   contractive Markov models or  autoregressive models with a Lipschitz mapping. For more thorough discussion of these examples, we  refer the reader to \cite{dedecker2004coupling, dedecker2005new} and Chapter 3 in \cite{dedecker2007weak}.  


As illustrated in the previous paragraph, the choice of the class ${\cal F}$ and  of the (semi-)norm  $ \|\cdot \|$ determine the type of dependence conditions for the  stochastic process. Following \cite{hang2017bernstein}, we consider the (semi-)norm 
\begin{align}
    \label{eq 2-7}
    ||f||=||f||_{\infty}+|f|_{SN} 
\end{align}
of a function $f:\mathbf{Z}\rightarrow\mathbb{R}$, where $\| f \|_\infty =\sup_{z\in\mathbf{Z}}|f(z)|$ and 
$|\cdot|_{SN}$ is a chosen semi-norm (for example,   $|f|_{\text{Lip}}$  or  $|f|_{TV}$). \cite{hang2017bernstein} made a literature review of some interesting examples of dynamic systems which are covered by the use of 
a semi-norms in Definition \ref{def 1}.

\section{Exponential Inequalities} \label{sec3}
\def\theequation{3.\arabic{equation}}	
   \setcounter{equation}{0}
   
In this section, we will derive several new  exponential inequalities for the mean $\frac{1}{n}\sum_{i=1}^{n}h(Z_i)$  
from a $\mathcal {C}_{\infty,\mathcal{F}}$ or $\mathcal{C}_{1,\mathcal{F}}$-mixing proces  $(Z_i)_{i \in \mathbb{N}}$, where $h: \mathbf{Z} \to \mathbb{R}$ is a function such that $E[h(Z_i)]=0$ ($i=1, \ldots , n$) and $h \in {\cal F}$ is allowed to depend on $n$. For $\mathcal{C}_{\infty,\mathcal{F}}$-mixing process, we show a slightly sharper inequality under more general conditions
compared to the pioneering work of \cite{hang2017bernstein}.  For (geometrical) $\mathcal{C}_{1,\mathcal{F}}$-mixing process, we show a sharp exponential inequality under a very mild assumption on the characteristic functions of $h(Z_i)$, which is much weaker than the Cram{\'e}r condition used in Edgeworth expansion. 
To our best  knowledge, this is the first  sharp inequality under such a general mixing condition. Furthermore, the  inequality holds for random variables taking values in an arbitrary measurable space $\mathbf{Z}$ and 
the value of $t$ in the event $\{|\frac{1}{n}\sum_{i=1}^{n}h(Z_i)|>t\}$ only appears in the exponential and not in
an additional remainder term. Note that this is desirable property that is often enjoyed only under the $\beta$-mixing (absolute regularity) condition and its special cases.

\subsection{$\mathcal{C}_{\infty,\mathcal{F}}$-mixing Processes}
\label{sec31}
For a semi-norm of the form \eqref{eq 2-7} Definition \ref{def 1}  with $p=\infty$  yields  a generic covariance inequality  for a 
$C_{\infty,\mathcal{F}}$-mixing process, that is   
\begin{align}
\label{eq 3.1-1}
    \max_{n}|\text{Cov}(Y,f(Z_{n+k}))|\leq C(k)||Y||_{1}(||f||_{\infty}+|f|_{SN}).
\end{align}

This coincides with the definition of $\mathcal{C}$-mixing process in \cite{hang2017bernstein} (see equation (2.7) therein), who derived
an exponential inequality 
for the mean $\frac{1}{n}\sum_{i=1}^{n} h(Z_i)  $
(where $\mathbb{E}[Z_i]  =0 $). Their work was motivated by the fact that the available results for the concentration of $ g(Z_1,..,Z_n)-E[g(Z_1,..,Z_n)]$
\citep[see, for example,][]{dedecker2005new,Chazottes_2005,Chazottes2012} do not use any variance information and can be improved if the variance of $h(Z_i)$ is small. However, as pointed out by \cite{hang2017bernstein},
such situations  appear quite frequently in  the analysis of learning algorithms. 
Another application, where  Bernstein-type inequalities reflect the effect of variance are useful,  are  uniform convergence rates for  kernel estimates 
in nonparametric regression
\citep[see][]{hansen2008uniform,kristensen2009uniform,li2012local,vogt2012nonparametric}.

 \cite{hang2017bernstein} obtained  a Bernstein-type inequality using the estimate \eqref{eq 3.1-1} iteratively to bound the Laplace transform of averaged partial sums after a “sub-lattice” trick. However, as we point out below, their approach yields an additional but unnecessary logarithmic factor and requires the stationarity of the process $(Z_{n})_{n\in\mathbb{N}}$. We obtain an improvement using a more direct  procedure to bound  Laplace transforms and  making full use of the fact that the upper bound of the covariance inequality \eqref{det1} only relies on $L_1$-norm of random variable $Y$.   Before we state our results, we introduce the following key assumption on the semi-norm $|\cdot|_{SN}$, which is also made in \cite{hang2017bernstein}.
 \begin{assumption}
  {\it \label{as 1}
 Given a semi-norm $|\cdot|_{SN}$ and a function $f $ with $ |f|_{SN}< \infty $, we assume 
    \begin{align*}
        |e^{f}|_{SN}\leq ||e^{f}||_{\infty}|f|_{SN}. 
    \end{align*}}
    \end{assumption}


\begin{example} ~
\label{example}
{\rm \cite{hang2017bernstein} listed some examples of semi-norms which 
satisfy Assumption \ref{as 1} and  
yield the $\mathcal{C}$-mixing ($\mathcal{C}_{\infty,\mathcal{F}}$-mixing) property of important  dynamical systems. These include  dynamic systems with piecewise-expanding maps, uniformly hyperbolic attractors and non-uniformly hyperbolic uni-modal maps
\citep[see also][]{viana1997stochastic}.
Theorem 2.9 in \cite{hang2017bernstein} establishes sufficient conditions for $\mathcal{C}$-mixing-type dynamic systems associated with these semi-norms and 
 for the sake of completeness we list these semi-norms here. 
 
\begin{itemize}
 \item[(1)] If $\mathbf{Z}$ is any given set we define the  semi-norm $|f|_{SN}=0$ for any function  $f:\mathbf{Z} \to \mathbb{R} $.
    \item[(2)] If $\mathbf{Z} \subset \mathbb{R}^{D}$ is open, we define for any  continuously differentiable 
     $f: \mathbf{Z} \to \mathbb{R} $  
     $$
     |f|_{SN}=\sum_{i=1}^{D}||\partial_{i}f||_{\infty},$$ where $\partial_{i}f$ denotes the partial derivative of $f$ with respect to its $i$-th variable.   
   \item[(3)] If $\mathbf{Z}=\mathbb{R}^{D}$ and  $\alpha\in (0,1]$, 
   we define for any 
   $\alpha$-Hölder continuous function 
    $f: \mathbf{Z} \to \mathbb{R} $ 
    \begin{align*}
        |f|_{SN}=\sup_{z\neq z'}\frac{|f(z)-f(z')|}{||z-z||_{E}^{\alpha}},
    \end{align*}
    where $||\cdot||_{E}$ denotes some norm on $\mathbb{R}^{D}$.
    \item[(4)]  If  $\mathbf{Z}=\mathbb{R}$, we define 
    for any function 
    $f: \mathbf{Z} \to \mathbb{R} $ with bounded total variation $|f|_{SN}=|f|_{TV}$
\end{itemize}
}
\end{example}

\begin{theo}
{\it 
\label{Th 1} 
Suppose that $(Z_{n})_{n\in\mathbb{N}}$ 
is a  $\mathbf{Z}$-valued  
  geometrically   $\mathcal{C}_{\infty,\mathcal{F}}$-mixing stochastic process  of order $(\nu, b, \gamma)$. Let  $ h \in {\cal F} $ satisfy Assumption \ref{as 1}, $E[h(Z_i)]=0$, $E(h^2(Z_i))\leq \sigma^{2}$ for $i=1, \ldots , n$,  $||h||_{\infty}\leq A$ and $|h|_{SN}\leq B$, for constants $\sigma, A ,B>0$. For (arbitrary) $\omega>1$ define $N_{0}:=\min\{n\in\mathbb{N}^{+}:\frac{A+B}{An^{\omega-1}}\leq \frac{1}{3}\}$, then for any $t>0$ and $n\geq N_0$ 
	\begin{equation}
		\label{eq 3.1-4}
		\mathbb{P}\Big \{ \Big |\frac{1}{n}\sum_{i=1}^{n}h(Z_{i})\Big |\geq t \Big\}\leq 2e^2\exp\Big (-\frac{nt^{2}}{4(\frac{\omega}{{b}}\log_{\nu}n)^{{1}/{\gamma}}(\sigma^{2}+tA)}\Big ) ~.
	\end{equation}}
\end{theo} 
\noindent 
It is of interest to compare this result with 
Theorem 3.1 in \cite{hang2017bernstein}.
First, note that Theorem \ref{Th 1} uses separate bounds $A$ and $B$ on the norm $\|\cdot \|_\infty $ and the semi-norm $|\cdot |_{SN}$, respectively,  and  that it is applicable for any $\omega >1$. Therefore, if $A$ is a fixed positive constant, 
 it gives an exponential bound,  whenever $|h|_{SN}$ is of order $O(n^{\alpha})$, for arbitrary fixed $\alpha>0$. In contrast, Theorem 3.1 in \cite{hang2017bernstein} provides such a bound only in the case  $\alpha\leq 2$. Second, Theorem \ref{Th 1} does not require the assumption of a stationary process 
 $(Z_{n})_{n\in\mathbb{N}}$. Third,
 the inequality \eqref{eq 3.1-4} is on the logarithmic level sharper compared to Theorem 3.1 in \cite{hang2017bernstein}. 
 \smallskip

 We now continue deriving a concentration inequality for algebraic $\mathcal{C}_{\infty,\mathcal{F}}$-mixing processes.     
\begin{theo} 
{\it \label{Th 2}
 Suppose $(Z_{n})_{n\in\mathbb{N}}$ is a $\mathbf{Z}$-valued algebraically $\mathcal{C}$-mixing process of order $\gamma$, where $\gamma>0$. For all $h:\mathbf{Z}\rightarrow\mathbb{R}$, we assume $E[h(Z_i)]=0$, $E[h^2(Z_i)] \leq \sigma^{2}$, for $i=1,...,n$, 
  $||h||_{\infty}\leq A$ and $||h||\leq B \leq n^{\alpha}$, for some $\sigma$, $A$, $B >0 $ 
 and  $\alpha>0$. 
If  Assumption \ref{as 1} is satisfied for $h$, we have for any $t>0$ and $n\geq N_{1}=\min\{n\in\mathbb{N}:n^{\alpha}>A\}$, 
\begin{equation}
    \label{eq 3.1-5}
    \mathbb{P}\Big \{ \Big  |\frac{1}{n}\sum_{i=1}^{n}h(Z_{i})\Big  |\geq t \Big \}\leq C_{1}\exp\Big (-\frac{n^{{\gamma}/{(\gamma+1})}t^2}{2(\sigma^{2}+tA)}\Big ), 
\end{equation} 
where $C_1=\exp \left(1+\frac{2e}{A}\right)$ .
}
\end{theo} 
 For algebraically ${\cal C}$-mixing processes there does not exist a directly comparable result. 
To our best knowledge, \cite{blanchard2019concentration} is  the only work  providing a concentration inequality for a sum  of 
Banach-valued  random variables from an  
algebraically  $\mathcal{C}$-mixing process. In contrast,  our Theorem \ref{Th 2} provides such an
inequality for the process $(h(Z_n))_{n\in \mathbb{N}}$ and $(Z_{n})_{n\in\mathbb{N}}$ is allowed to take values in any given measurable space $\mathbf{Z}$.
Although both papers consider different scenarios,  the derived 
bound in \eqref{eq 3.1-5} 
shows similarities with  the bound  obtained by   \cite{blanchard2019concentration} in their Theorem 3.5 and Proposition 3.6. In principle,  their results can be applied 
to  processes of the form $(h(Z_{n}))_{n\in\mathbb{N}}$. However, 
in order to do this, one has to show that the $(h(Z_{n}))_{n\in\mathbb{N}}$ remains  a  $\mathcal{C}$-mixing process and its $\mathcal{C}$-mixing coefficient is smaller than that of $(Z_{n})_{n\in\mathbb{N}}$. To our best knowledge, there are no general sufficient conditions that guarantee this property. 
Showing such results seems to be inherently difficult as indicated by \cite{dedecker2007weak} for $\tilde{\alpha}$-mixing processes.
Up to our knowledge, the most recent progress is that, for a real valued 
$\tilde{\alpha}$-mixing process
$(Z_{n})_{n\in\mathbb{N}}$, the process $(h(Z_{n}))_{n\in\mathbb{N}}$ is still  $\tilde{\alpha}$-mixing  with not enlarged mixing coefficient, if the function $h$ is strictly monotone \citep[see Chapter 3 in][]{dedecker2007weak}.


\subsection{$\mathcal{C}_{1,\mathcal{F}}$-mixing processes }
\label{sec32}
\par In this section we will derive a sharp exponential inequality reflecting the variance information for a  geometrically   $\mathcal{C}_{1,\mathcal{F}}$-mixing $\mathbf{Z}$-valued  process $(Z_n)_{n\in\mathbb{N}}$.  In this case, we have  $p=1$ and   with \eqref{eq 2-7}  Definition 
\ref{def 1} yields the
 covariance inequality 
\begin{align}
    \label{eq 3.2-1}
    \max_{n\in\mathbb{N}}\text{Cov}(Y,f(Z_{n+k}))\leq C(k)||Y||_{\infty}(||f||_{\infty}+|f|_{SN}).
\end{align}
Compared  with \eqref{eq 3.1-1},  this inequality is much weaker since only the $\|\cdot \|_\infty$-norm of the random variable $Y$ is involved. This simple difference 
makes the derivation of  exponential inequalities for ${1 \over n} \sum_{i=1}^n h(Z_i) $ from a (general) $\mathbf{Z}$-valued 
$\mathcal{C}_{1,\mathcal{F}}$-mixing process $(Z_i)_{i \in \mathbb{N}}$ much harder:  
 First, it is nearly impossible to  use similar arguments as given  in the proof of Theorems \ref{Th 1} and \ref{Th 2}, since 
the use of the $L_1$-norm of $Y$ in the application of the covariance inequality is crucial  for deriving bounds on the Laplace transform.
Second, the other powerful tool in showing sharp exponential inequalities for weakly dependent processes is coupling.  While, there 
are many strong coupling lemmas available  for processes with mixing coefficients defined for  two  sigma fields \citep[see][for  a nice survey] {merlevede2002coupling}, 
the literature on coupling for processes with 
mixing coefficients defined for a sigma field and a random variable is rather scarce. We refer to \cite{rio1995functional}, \cite{rio2017asymptotic}, \cite{peligrad2002some} and \cite{dedecker2005new}, who derived  some result for  real valued $\alpha$, $\Tilde{\alpha}$- and $\tau$-mixing  processes.  Some extensions for  $\tau$-mixing process taking values in a Polish space can be found \cite{dedecker2006parametrized}. However, in this paper, we investigate  $\mathbf{Z}$-valued $\mathcal{C}_{1,\mathcal{F}}$-mixing processes, where $\mathbf{Z}$ is taking values in  a general metric space. Thus, the model under consideration  is substantially more general (with respect to the image set of the random variables and the dependence concept) as in the cited references. Third, as mentioned before, coupling of $h(Z_{n})$ does not always work well, because, to the  best of our knowledge, it is not obvious if the mixing coefficients 
corresponding to the pair $(\mathcal{M},h(Z))$ can be bounded by that corresponding   to $(\mathcal{M},Z)$. \\
Summarizing, the available tools can not be used to derive exponential inequalities valid for the whole class of $\mathcal{C}_{1,\mathcal{F}}$-mixing processes without any additional assumptions on measurable space $\mathbf{Z}$. 
However, we will show below that a mild assumption (Assumption \ref{as 3}) on the characteristic function of $h(Z_i)$ is sufficient to bypass this barrier.

\begin{assumption}
\label{as 2}
{\it For the semi-norm $|\cdot|_{SN}$ and function $f \in {\cal F}$, we assume  
    \begin{align}
    \label{eq 3.2-2}
       S(u) := |\cos(uf)|_{SN}\lor |\sin(uf)|_{SN}\leq C_{B} |u|^{\eta_{0}}|f|_{SN},
    \end{align}
    for  every $u\in\mathbb{R}$, 
   where $\eta_{0} \geq  0$ and $C_{B} > 0 $ are constants independent of sample size $n$.}
\end{assumption}

By straightforward calculations we can show the semi-norms listed in Example \ref{example} satisfy Assumption \ref{as 2}.
 Note  that the semi-norm $|\cdot|_{SN}$ is  determined by the particular   $\mathcal{C}_{1,\mathcal{F}}$-mixing condition  under consideration and  often chosen to  reflect smoothness. 
 Roughly speaking, Assumption \ref{as 2} 
    essentially requires that the function $S(u)$ diverges at a polynomial rate as $|u|\nearrow \infty$. Another key point is that it does not impose any restriction on the boundedness of the function $f$, which is much more general than Assumption \ref{as 1}.
We finally introduce the following assumption about the characteristic functions
   $\psi_{h(Z_i)} (u)  = E [ \exp ( i u h(Z_i) )] $  of the random variables $h(Z_i)$ ($i=1, \ldots , n$). 
\begin{assumption}
\label{as 3}{
\it  Consider a
   $\mathbf{Z}$-valued  
  geometrically   $\mathcal{C}_{1,\mathcal{F}}$-mixing  process  $(Z_{n})_{n\in\mathbb{N}}$. 
We assume there exists a subset $I_n'\subset I_n:=\{1,2,\ldots ,n\}$ such that the following two conditions are satisfied.
\begin{itemize}
    \item [1)] 
    There exists constants $0<c<C\leq 1$ such that, for every given $n\in\mathbb
    N$, $\frac{\text{Card}(I_n')}{n}\in [c,C]$.
    \item [2)] For each $i\in I_n'$ and any fixed $\kappa, \Theta>0$ that are independent of $n$,
    there exists a constant 
    $N_0\in\mathbb{N}$ and a sequence $\{\tau_n \}_{n \geq N_0}$ independent of $i$ such that
    $$
    \tau_{n}\geq \frac{2\kappa(\log n)(b^{-1}\Theta\log n)^{\frac{1}{\gamma}}}{n},
    $$ 
    (here  $\gamma$ is is the constant in the definition of the  geometrically   mixing process)
    such that, for any $n\geq N_0$,  
    \begin{align}
    \label{cramer}
        |\psi_{h(Z_i)}(u)|\leq 1-\tau_{n},\ \forall\ |u|\in [M_n,\infty),
    \end{align}
     where $\{M_n\}_{n \geq N_0} $ is a sequence of increasing positive constants satisfying $M_{n}\leq n^{\beta}$, for some $\beta>0$.
\end{itemize} 
If $(Z_{n})_{n\in\mathbb{N}}$ is identically distributed, we only require point 2)  above.
}
\end{assumption}
Assumption \ref{as 3} is weaker than the famous Cram{\'e}r condition 
 $$
{\limsup}_{|t|\rightarrow \infty}|\psi_{h(Z_i)}(t)|<1,
 $$
 which is frequently used in the field of Edgeworth expansions 
\citep[see, for example,][]{hall2013bootstrap}.
  Note that it is is implied by  the weak Cram{\'e}r condition introduced by \cite{angst2017weak}. Thus, all the examples provided in this references  do not satisfy    Cram{\'e}r's  condition but  satisfy \eqref{cramer}. Moreover, for statistical applications, particularly kernel smoothers, we mention the recent work of \cite{calonico2018effect} who introduced an $n$-dependent Cram{\'e}r condition (see Assumptions S.I.3.3 and S.II.3.3 in the corresponding online supplement), which implies in fact assumption \eqref{cramer}.  We further show that, under very general  conditions on the kernel, the summands in common kernel density estimator satisfy  condition \eqref{cramer}; see Proposition \ref{prop A2} below. This allows us to use the following Theorem \ref{thm 3} to reveal the concentration phenomenon of kernel density estimators under a  geometrically $\mathcal{C}_{1,\mathcal{F}}$-mixing  assumption, which is usually a  crucial step in obtaining the convergence rates of many density-based plug-in type estimators, like  the local conditional mode set estimator discussed in Section \ref{sec4}. Similar results as in  Proposition \ref{prop A2} can be derived for 
  the summands   in kernel regression estimators, but the details are omitted for the sake of brevity.

\par  For the statement of the main theorem of this section,  we define the decomposition of index set $\{ 1, \ldots , n \} = \cup_{j=1}^{p_n}I_j$ with $p_n \le n/2$ by
\begin{align}
\label{det6}
 I_{j} & =\left\{
\begin{aligned}
&\{j+kp_{n}:0\leq k\leq [n/p_n]\}  , & \text{when $1\leq j\leq r_{n}$}, \\
&\{j+kp_{n}:0\leq k\leq [n/p_n]-1\} , & \text{when $r_n +1\leq j\leq p_{n}$} ,
\end{aligned}
\right. 
\end{align}
where $r_{n}=n-p_n \lfloor n/p_n \rfloor$.

\begin{theo} \label{thm 3}
\it 
Suppose that $(Z_{n})_{n\in\mathbb{N}}$ 
is a  $\mathbf{Z}$-valued  
  geometrically   $\mathcal{C}_{1,\mathcal{F}}$-mixing  process  of order $(\nu, b, \gamma)$. Let  $ h \in {\cal F} $ such that for $i=1, \ldots , n$  the random variable $h(Z_i)$   is centered and  satisfies a
Bernstein condition, that is $E|h(Z_{i})|^{k}\leq \frac{1}{2}k!A^{k-2}E|h(Z_{i})|^2$ for all $k\geq 3$
(with a universal constant $A>0$).
Further assume that  $|h|_{SN}\leq n^{B}$, for some $B>0$.
Then,  under Assumptions \ref{as 2} and \ref{as 3}, there exists an integer $N_0 \in \mathbb{N}$ such that, for any $n\geq N_0$, $t>0$ and $\delta>1$,
\begin{align}
\label{eq 3.2-3}
    \mathbb{P}\Big ( \Big |\frac{1}{n}\sum_{i=1}^{n}h(Z_{i}) \Big |>t \Big )\leq 2p_n\exp\Big  (-\frac{nt^2}{p_n(\sigma^{2}_{\max}+{At}/{3})}\Big  )+\frac{C'_{\delta} + 8\max_{i=1}^n E|h(Z_i)|}{n^{\delta-1}},
\end{align}
where $p_{n}=(b^{-1}\Theta_\delta \log_{\nu}n)^{\frac{1}{\gamma}}$, 
$$
\sigma_{\max}^{2}=\max_{1\leq j\leq p_n}\frac{1}{|I_j]}\sum_{k=1}^{|I_j]}{\rm Var}\left(h(Z_{j+kp_{n}})\right)
$$
and
$\Theta_\delta$ and
$C'_{\delta}$ are universal constants specified in equation \eqref{1} and \eqref{2} in the appendix, respectively.  Furthermore, if  $Z_{1}, \ldots , Z_n$ are identically distributed, the inequality \eqref{eq 3.2-3} simplifies to 
\begin{align}
    \label{eq 3.2-4}
     \mathbb{P}\Big  (\Big  |\frac{1}{n}\sum_{i=1}^{n}h(Z_{i})\Big  |>t\Big  )\leq 2p_n\exp\Big  (-\frac{nt^2}{p_n({\rm Var}(h(Z_1))+{At}/{3})}\Big )+\frac{C'_{\delta} + 8E|h(Z_1)|}{n^{\delta-1}}.
\end{align}
\end{theo}

We emphasize that the inequalities   \eqref{eq 3.2-3} and \eqref{eq 3.2-4} hold for any $\delta >0$.
This constant only affects the sharpness of the exponential part on the right-hand side  up to a constant factor (entering through $p_n$). On the other hand, increasing 
the constant $\delta$ makes the  second term smaller  (but the constants $C'_{\delta}$ and $\Theta_{\delta}$ increase with $\delta$). 
Note also that this remainder is independent of $t$. This makes the estimates \eqref{eq 3.2-3} and \eqref{eq 3.2-4} different from many exponential inequalities for the mean of dependent random variables, where the remainder term  does in fact depend on  $t$ \citep[see, for example Theorems 1.3 and 1.4 in][]{bosq2012nonparametric}. As we pointed out before, this is a property often enjoyed by $\beta$-mixing process and its spacial cases only.  We also emphasize that the  assumptions in Theorem \ref{thm 3} on the dependency structure are very general and that our result is valid for any set $\mathbf{Z}$.

\section{ Conditional Mode Set Estimation }
\label{sec4}
\def\theequation{4.\arabic{equation}}	
   \setcounter{equation}{0}
   
In this section we will use the new exponential inequalities to derive convergence rates for  a kernel-based estimator of the local conditional mode set of the density of a stochastic process satisfying a geometrically  $\mathcal{C}_{1,\mathcal{F}}$-mixing condition. The conditional mode is an alternative  regression method to explore the relationship between an output $Y$ and an input $X$. As pointed out by \cite{chen2016nonparametric}, the conditional mode  often yields a better  reflection of the trend of the output variables and a  sharper prediction band. Moreover, it is also more robust with respect to the skewness of the distribution of the error. 
\smallskip 

Suppose $(X,Y)$ is a random vector taking values in set $\mathbb{X}\times \mathbb{Y} \subset \mathbb{R}^{D_X} \times  \mathbb{R}^{D_Y}$ with joint density $f_{(X,Y)}$. Suppose that the  marginal density $f_{X}$  of  $X$ satisfies 
\begin{align}
\label{det7}
\inf_{x\in\mathbb{R}^{D_X}}f_X(x)>0,
\end{align}
and define 
\begin{align}
\label{eq 4-1}
    \Omega_{x_0}=\{y\in\mathbb{Y}: (x_0,y)\in\mathbb{X}\times \mathbb{Y}\} \subset \mathbb{R}^{D_Y}. 
\end{align}
For any  $x_0 \in \mathbb{X}$, we define the “\textit{local conditional mode set at point $x_0$}” by
\begin{align*}
    \mathcal{Y}_{x_0}&= \Big \{y\in\mathbb{Y}:\exists\ \text{cube}\ I_y\subset \Omega_{x_0}\ \text{such that}\  f_{Y|X}(y|x_0)=\arg\max_{z\in I_y}f_{Y|X}(z|x_0) \Big\},\\
    &=\Big\{y\in\mathbb{Y}:\exists\ \text{cube}\ I_y\subset \Omega_{x_0}\ \text{such that}\  f_{(X,Y)}(x_0,y)=\arg\max_{z\in I_y}f_{(X,Y)}(x_0,z)\Big \},
\end{align*}
where $f_{Y|X}$ denotes the conditional density of $Y$ given  $X=x_0$ and the equality between holds because of assumption  \eqref{det7}. Each element of the set $\mathcal{Y}_{x_0}$ is a local maximizer of the function $f_{(X,Y)}(x_0,z)$. We introduce the following key assumption on the mode set. 
\begin{assumption}
    {\label{assumption 4-1}
   \it For each given $x_0\in\mathbb{X}$, 
   the local conditional mode set $\mathcal{Y}_{x_0}$ is a finite set with $K_{x_0}$ elements. We further assume that 
   $\mathbf{K}= \sup_{x_0\in \mathbb{X}}K_{x_0} < \infty  $.}
  \end{assumption}
\par Note that in contrast to most of the literature we do not assume that the local conditional mode is unique. A uniqueness assumption makes the local conditional mode set equal to global mode set and simplifies the analysis of conditional mode estimators substantially. 
Under the uniqueness assumption, \cite{collomb1986note} proposed to estimate the global mode
directly by determining the maximum of  a kernel density estimator. Moreover, this assumption implies that the global conditional mode, 
\begin{align}
\label{eq 4-6a}
\theta(x):=\text{Mode}(Y|X=x_0)=\arg\max_{y}f_{Y|X}(y|x_0), 
\end{align}
is a function of $x$ and one can develop a framework similar to conditional mean or quantile regression considering the model 
\begin{align}
\nonumber
    &Y=\theta(X)+\epsilon,
\end{align}
where the error satisfies $\text{Mode}(\epsilon|X=x)=0.$
This setting is obviously more robust to the skewness and outliers and meanwhile several authors have 
worked on estimating the function $\theta$ under different model assumptions using  traditional  regression techniques \citep[see][among others]{lee1989mode, lee1993quadratic,kemp2012regression,yao2014new,krief2017semi,ota2019quantile,zhou2019bandwidth,feng2020statistical,kemp2020dynamic,ullah2021modal,ullah2022nonlinear,ullah2023semiparametric}. However, in many applications the uniqueness assumption ($\mathbf{K}=1$) might be too restrictive. If $\mathbf{K}>1$, the function $\theta$ in \eqref{eq 4-6a} is  a set-valued mapping. Compared to the literature on mode estimation  under the uniqueness assumption, the work  considering the case $\mathbf{K}>1 $ is much more scarce. To our best knowledge, the only available theoretical work in this direction is \cite{chen2016nonparametric} who investigated the asymptotic properties of an estimator of the local conditional mode set under the assumption of the existence of the second-order partial derivatives of the joint density of $(X,Y)$. As in \cite{collomb1984proprietes}, their 
 set estimator is based on kernel estimators of the joint density  and its higher-order partial derivatives.

Nearly all of the cited literature investigates estimators under the assumption of iid data. However, in many applications, such an assumption might not be reasonable. For example, one might be interested in the most likely value of $W_{n+1}$ based on the observation of $W_n$ from a stationary time series $(W_n)_{n\in \mathbb{N}}$. In this case we are interested in the mode set of the conditional density of $Y_i=W_{i+1}$ given  $X_i =W_i$.
Therefore, in this section, we  investigate the convergence rate of the Hausdorff loss of a kernel-based local conditional mode set estimator for dependent data. 
Moreover, in contrast to \cite{chen2016nonparametric}, 
we do not assume that the conditional density function is  differentiable at its local or global maxima and our approach therefore avoids estimating the first and second order partial derivatives of the conditional density. Furthermore, we also consider the case of vector-valued output $Y$.
\par To be precise, consider a $\mathbb{X}\times \mathbb{Y} \subset \mathbb{R}^{D_X} \times  \mathbb{R}^{D_Y}$-valued stochastic process $(X_i,Y_i)_{i \in \mathbb{N}}$ with identical marginal distributions with  density $f_{(X;Y)}$.  Let $K$ denote a (one-dimensional) kernel function supported on the interval $[-1,1]$, and define the product kernels
\begin{align} \label{det21}
 K_{x_0,h}(X_i)&=\prod_{s=1}^{D_X}K\Big(\frac{X_{is}-x_{0s}}{h}\Big),\ K_{y,h}(Y_i)=\prod_{s=1}^{D_Y}K\Big(\frac{Y_{is}-y_{s}}{h}\Big).
\end{align}
where $x_0=(x_{01},...,x_{0D_X})^\top $, 
$X_i=(X_{i1},...,X_{iD_X})^\top \in\mathbb{R}^{D_X}$, $Y_i=(Y_{i1},...,Y_{iD_{Y}})^\top\in\mathbb{R}^{D_Y}$ and $h>0$ is a bandwidth. We consider the common product-type estimator 
\begin{align} \label{det19}
    &\hat{f}_{(X,Y)}(x_0,y)=\frac{1}{nh^{D}}\sum_{i=1}^{n}K_{x_0,h}(X_i)K_{y,h}(Y_i),
\end{align}
for the joint density of the vector $(X_i,Y_i)$, where $D=D_X+D_Y$ and define
    \begin{align}
    \label{eq 4-8}
    \widehat{\mathcal{Y}}_{x_0}=\{y\in\Omega_{x_0}:\ \exists\ \text{cube}\ I_y\subset \Omega_{x_0}\ \text{such that}\  \hat{f}_{(X,Y)}(x_0,y)=\arg\max_{z\in I_y}\hat{f}_{(X,Y)}(x_0,z) \},
\end{align}
as the estimator of the set of local conditional modes.

Similar to \cite{chen2016nonparametric}, we 
measure the performance of the estimator  $\widehat{\mathcal{Y}}_{x_0}$ by 
the Hausdorff loss  
\begin{align}
    \label{eq 4-9}
    \triangle_{n}(x_0):=\max
    \Big \{\sup_{y\in\widehat{\mathcal{Y}}_{x_0}}d(y,\mathcal{Y}_{x_0}), \sup_{y\in\ \mathcal{Y}_{x_0}}d(y,\widehat{\mathcal{Y}}_{x_0}) \Big \},
\end{align}
 where 
 $d(y,A)=\inf_{x\in A } ||y-x||_{E}$ denotes the distance of $y$ from the set $A$ and $||\cdot||_E$ is the Euclidean norm on $\mathbb{R}^{D_Y}$. To study $ \triangle_{n}(x_0)$ we make the following
regularity assumptions on the joint density $f_{(X,Y)}$.

 \begin{assumption} 
{\it 
 \label{as 5} 
   The joint density  
   $f_{(X,Y)}$ satisfies the following conditions 
    \begin{itemize}
        \item [(C1)] The joint density $f_{(X,Y)}(x,y)$ is bounded on $\mathbb{X}\times\mathbb{Y}$  and coordinate-wise  in the $\alpha$-Hölder class.  That is, if $f_j$  is the function obtained by considering $f$ as function of the $j$ coordinate and keeping all other coordinates fixed, then $f_j$  is $\lfloor \alpha \rfloor $  times differentiable and the $\lfloor \alpha \rfloor$-th derivative is Hölder continuous of order $(\alpha - \lfloor \alpha \rfloor)$; ($j=1, \ldots , D$); 
         see  Definition 1.2 in  \cite{tsybakov2004introduction}.
        \item [(C2)] 
     For any  $x_0 \in \mathbb{X} $, the function  $y \to  f_{(X,Y)}(x_0,y)$ is {$\beta$-distinguishable on its mode set for some $\beta >0$}, that is: for any $y^*\in\mathcal{Y}_{x_0}$, $\beta $ is the smallest positive number such that
        \begin{align*}
         \lim_{\delta\rightarrow0^+}  \inf_{
         \|y^*-y\|_{E} < \delta }  \frac{|f(x_0,y^*)-f(x_0,y)|}{\|y^*-y\|_{E}^{\beta}}\geq c_{y^{*}}>0,
        \end{align*}
        holds for some positive  constant $c_{y^*}$ that does not depend  of $y$.
    \end{itemize}
    }
     \end{assumption}
 \begin{rem}
 \label{remark}
 {\rm 
The constant  $\beta$  in condition (C2) geometrically describes  how “flat” the function is at its conditional mode point.  If  $\mathbb{Y}\subset \mathbb{R}$ and the function $y \to f_{(X,Y)}(x_0,y)$ 
is $l \geq 1$ times  continuously differentiable at $y^*$, the magnitude of $\beta$ is determined by the rate of $\big |\frac{\partial f_{(X,Y)}}{\partial y}(x_0,y) \big |$ as $y\rightarrow y^*$.
Furthermore, to illustrate  this condition in the case  $\mathbb{Y}\subset \mathbb{R}$, we list 3 specific examples.
\begin{itemize}
    \item[(1)] If $f_{(X,Y)}(x_0,y)$ is a (mixture) triangular-type or (mixture) exponential-type distributions, $f(x_0,y)$ is  distinguishable with $\beta=1$ on its mode set.
    \item[(2)]  
   Let  $X$  be uniformly distributed on the interval $[0.4,0.9]$ and $w_1, \ldots , w_K$
    are positive weights with $\sum_{k=1}^{K}w_k=1$ for some $w_k\geq 0$ and define 
    \begin{align}
    \label{eq 4-10}
        f_{Y|X}(y|x_0)=\sum_{k=0}^{K}\frac{w_{k}}{2}(1-|y-2k|^{x_0})1[|y-2k|\leq 1].
    \end{align}
  Then $f_{(X,Y)}(x_0,y)$ is distinguishable with   the $\beta=x_0 $ ($x_0\in [0.4,0.9]$) on its mode set $\{2k:k=0,1,..,K\}$.
    \item[(3)]  If $f_{(X,Y)}(x_0,y)$ is a (mixture) of Gaussian distributions (with respect to $y$), then $f_{(X,Y)}(x_0,y)$ is distinguishable with $\beta=2$ on its mode set.
\end{itemize}
}
 \end{rem} 
 We further introduce the following crucial assumptions of the kernel functions used in our mode estimator.
\begin{assumption}
    \label{kernel}
     {\it We assume the kernel $K$ in \eqref{det21} satisfies the following conditions.
     \begin{itemize}
         \item [(K1)] $K$ is a bounded, non-negative and symmetric density  with compact support $[-1,1]$.Further, $K$ is strictly decreasing on $[0,1]$ and, for any $\epsilon_n\searrow 0$, there exists a constant $\beta_1>0$ such that $K(1-\epsilon_n)\geq \epsilon_n^{\beta_1}$.  We also assume $K''$ is bounded.
         \item [(K2)] We assume $K$ is twice differentiable with  first and second order derivatives as $K'$ and $K''$, respectively.
         \item [(K3)] If $K'$ is not a negative constant on the interval $[0,1]$, we assume the $|K'|$ is non-decreasing on $[0,1]$ and $|K'(\epsilon_n)|=O(\epsilon_n^{\beta_2})$, for some $\beta_2>0$.
     \end{itemize}  }
\end{assumption}

 Many widely chosen kernels such as the triangle, Epanechnikov-, Silverman- and truncated Gaussian kernels satisfy this assumption. Moreover, all kernels of the form $(1-|u|^a)^b1[|u|\leq 1]$, $a,b\geq 1$, satisfy Assumption \ref{kernel}. Our first result indicates that, together with condition (C1) in Assumption \ref{as 5}, Assumption \ref{kernel} implies that the $K_{x_0,h}(X_i)K_{y,h}(Y_i)$ satisfies  condition \eqref{cramer}, which is crucial for the application of Theorem \ref{thm 3}. 
 

\begin{prop}
\label{prop A2}
  {\it Suppose that  $(Z_n)_{n\in \mathbb{N}} = ( (X_i,Y_i))_{n\in \mathbb{N}}$ are identically distributed random variables taking value in $\mathbb{X}\times\mathbb{Y}\subset \mathbb{R}^D$ with joint density $f_{(X,Y)}(x,y)$. If  $h=O(n^{-\eta})$ for some $\eta>0$ and  Assumption \ref{kernel} and condition (C1) of Assumption \ref{as 5} hold, the product kernel $K_{x,h}(X_1)K_{y,h}(Y_1)$ satisfies Assumption \ref{as 3}, for any $(x,y)\in\mathbb{X}\times\mathbb{Y}$.}
\end{prop}
Finally, we also need to specify Assumption \ref{as 2} for the application of  Theorem \ref{thm 3}. For this purpose, we introduce the following mild assumption.
\begin{assumption}
    \label{kernel 2}
 {\it The kernel  $K_{x_0,h}K_{y,h}$ defined by the product of the kernels in \eqref{det21} satisfies  
    $|K_{x_0,h}K_{y,h}|_{SN}\leq C_{D}n^{B}$,
where the constants $C_D$ and $B$ do not depend on $n$. }
 %
\end{assumption}
As pointed out in \cite{hang2018kernel}, Assumption \ref{kernel 2} holds for $\mathcal{C}$-mixing ($\mathcal{C}_{\infty,\mathcal{F}}$-mixing) processes. For example, when $|\cdot|_{SN}$ is chosen as Lipschitz norm, we can show the growth condition of the Lipschitz norm of the product of kernels is of order $O(\frac{1}{h})$.  Therefore, Assumptions \ref{kernel}, \ref{kernel 2} and condition (C1) of Assumption \ref{as 5} ensure that Theorem \ref{thm 3} can be applied to investigate the theoretical guarantee of the local conditional mode set estimator.
The following theorem is the main result of this section and gives  the convergence rate of the  mode set estimator \eqref{eq 4-8} with respect to the Hausdorff loss.

\begin{theo}
    {\it
\label{Th 4}
 Suppose $(X_i,Y_i)_{i\in \mathbb{N}}$ is a geometrically   $\mathcal{C}_{1,\mathcal{F}}$-mixing process of order $(\nu,b,1)$ with identically $1$-dimensional marginal distributions
 taking values  
 in set $\mathbb{X}\times \mathbb{Y} \subset \mathbb{R}^{D_X} \times  \mathbb{R}^{D_Y}$. Let $D=D_{X}+D_{Y}$ and
 $h=(\frac{\log n}{n})^{\frac{1}{2\alpha+D}}$. Then, under Assumptions \ref{as 5}, \ref{kernel} and \ref{kernel 2}, the Hausdorff loss of the mode estimator \eqref{eq 4-8} satisfies  
    \begin{align}
    \label{det2}
        \triangle_{n}(x_0)=\mathcal{O}_{\mathbb{P}}\Big ( \Big(\frac{\log n}{n}\Big)^{\frac{\alpha}{\beta(2\alpha+D)}}\Big ).
    \end{align}
   }
\end{theo}

Theorem \ref{Th 4} gives a precise statement how the  constant $\beta $
affects the convergence rate of the  Hausdorff loss $\triangle_{n}(x_0)$.
As we pointed out in Remark \ref{remark}, $\beta$ is related  the local smoothness of the  conditional density at its local mode. A smaller value of $\beta$ indicates more roughness and  leads to a faster rate in \eqref{det2}. Thus, unlike to the estimator of the conditional density at the mode, where more smoothness improves the rate of convergence, the  local mode set estimator performs better for smaller values of $\beta$. More specifically, when $\beta<1$, the convergence rate of $\triangle_{n}(x_0)$ is faster than the optimal  convergence rates of the  density estimator for $f_{(X,Y)}$.  Moreover, in the case $1\leq \beta<2 $, the statistic \eqref{eq 4-8} provides a consistent estimator of the local conditional mode set, while  the method of \cite{chen2016nonparametric}  is not  applicable.

More specifically, these authors proposed  an alternative estimate of the local conditional modes set in the case of iid data. Under the assumption of a four times continuously differentiable density $f_{X,Y)}$ and the condition
$|\frac{\partial_{}f_{(X,Y)}(x_0,y^*)}{\partial y^2}|>0$
they proved the rate
\begin{align}
   \label{det3} 
\mathcal{O}_\mathbb{P}\Big ( h^2 + \sqrt{1 \over n h^{D+2}}
\Big ) = \mathcal{O}_\mathbb{P} \Big ( \Big ( { 1 \over  n} \Big ) ^{2/(7+D_X)}
\Big ) 
\end{align}
where the last equality holds for the bandwidth $h=n^{-1/(7+D)}$.  The assumptions in \cite{chen2016nonparametric} correspond to $\alpha =4 $
and $\beta = 2 $ in condition (C1) and (C2), and, by Theorem \ref{Th 4} we obtain the rate 
$$
\mathcal{O}_{\mathbb{P}}\Big (\Big(\frac{\log n}{n}\Big)^{{2}/{(9 + D_X)}}\Big )
$$
for the estimator proposed in this paper, 
which is (slightly) worse  than the rate \eqref{det3} for the estimator of \cite{chen2016nonparametric}. 
In this sense, Theorem \ref{Th 4}
provides statistical guarantees 
for an alternative estimator  of  the  set of local conditional modes from dependent data. In some cases this estimate 
also improves the current state of the art \citep{chen2016nonparametric} for iid data.

\section{Conclusions and Discussion}
In this paper, we derived sharp exponential inequalities for an empirical-mean-type statistics of the form \(\frac{1}{n} \sum_{i=1}^{n} h(Z_{i})\),
where \( (Z_{i})_{i\in\mathbb{N}} \)  taking values in a measurable space \( \mathbf{Z}^{\mathbb{N}} \) 
with a very general dependence structure. We introduced the 
concept of \( \mathcal{C}_{p,\mathcal{F}} \)-mixing  $(1 \leq p \leq \infty$)
as  a general dependence condition that includes most of the currently considered non-strong-mixing processes and dynamical systems as special cases. While a lower value for \( p \) implies greater generality (which means a weaker dependence condition), the most restrictive condition is obtained for the choice $p= \infty $, which yields  \( \mathcal{C}_{\infty,\mathcal{F}} \)-mixing. This dependence concept coincides  with \( \mathcal{C} \)-mixing  which has  been discussed before in the literature. 
\\
The exponential inequality for  \( p = \infty \)  refines previous pioneering results in \cite{hang2017bernstein} by providing a sharper estimate under more general conditions. For the case  \( p = 1 \)
we provide a sharp exponential inequality 
under a very mild assumption on the characteristic function of \( h(Z_i) \)  for all geometrically \( \mathcal{C}_{1,\mathcal{F}} \)-mixing processes without requiring specific assumptions on the set \( \mathbf{Z} \).
To the best of our knowledge, this result is  the first in this direction.  
As an illustration   we use this inequality to investigate the convergence rate of a local conditional mode set estimator with respect to Hausdorff loss when the data is coming from a geometrically  \( \mathcal{C}_{1,\mathcal{F}} \)-mixing process and the marginal distribution does not necessarily have a differentiable density.

\acks
This work was supported by TRR 391 \textit{Spatio-temporal Statistics for the Transition of Energy and Transport} (Project number 520388526) funded by the Deutsche Forschungsgemeinschaft (DFG, German Research Foundation).

\newpage

\appendix
\section{Proof of Main Results}
\renewcommand{\theequation}{\thesection.\arabic{equation}}
\setcounter{equation}{0}

Before we start to prove  Theorem \ref{Th 1}, \ref{Th 2} and  \ref{thm 3}, we first introduce the following sub-lattice method which will  be used in the proof of all theorems in Section \ref{sec3}. Let $p_{n}$ be a positive integer no larger than $n/2$. By denoting $l_{n}=[n/p_{n}]$ and $r_n=n-l_n p_n$, for $j=1,2,\dots,p_n$, we define
\begin{equation}
\label{eq a-1}
I_{j}=\left\{
\begin{aligned}
&\{j+kp_{n}:0\leq k\leq l_{n}\}  , & \text{when $1\leq j\leq r_{n}$}, \\
&\{j+kp_{n}:0\leq k\leq l_{n}-1\} , & \text{when $r_n +1\leq j\leq p_{n}$}.
\end{aligned}
\right.
\end{equation}
Then, we have
\begin{align}
    \label{eq a-2}
    \frac{1}{n}\sum_{i=1}^{n}h(Z_{i})=\sum_{j=1}^{p_{n}}w_{j}\sum_{k=0}^{|I_j|-1}\frac{h(Z_{j+kp_n})}{|I_j|},
\end{align}
where
$$
\ w_{j}=\frac{|I_j|}{n}
$$
and  $\sum_{j=1}^{n}w_{j}=1$. Now we start to prove the results in Section 3.\\\\

{\bf Proof of Theorem \ref{Th 1}.} According \eqref{eq a-2}, we have
\begin{align}
\label{eq a-3}      &E\left[\exp\left(\frac{\lambda}{n}\sum_{i=1}^{n}h(Z_i)\right)\right]=E\left[\exp\left(\lambda\sum_{j=1}^{p_n}w_j\sum_{k=0}^{|I_j|-1}\frac{h(Z_{j+kp_n})}{|I_j|}\right)\right]
\leq \sum_{j=1}^{p_n}w_j E_j
\end{align}
where 
\begin{align}
    \label{hd1}
    E_j  & =E\left[\exp\left(\lambda\sum_{k=0}^{|I_j|-1}\frac{h(Z_{j+kp_n})}{|I_j|}\right)\right].
\end{align}
Therefore we  bound the left hand side of \eqref{eq a-3} by a deriving  a (uniform) upper bound for $E_1, \ldots , E_{p_n}$. 
For this purpose, we define, for  $j\in\{1,2,\dots,p_n\}$ and $l\in\{0,1,2,\dots, |I_j|- 1\}$,
\begin{align} 
\label{hd2} 
g_j (z) & =\exp\left(\frac{\lambda h(z)}{|I_j|}\right),  \\
\label{hd3} 
S_{l}& =\exp\left(\lambda\sum_{k=0}^{l}\frac{h(Z_{j+kp_n})}{|I_j|}\right).
\end{align}
Since $||h||_{\infty}\leq A$, it follows that $||S_l||_{\infty}\leq \exp(\lambda lA/|I_j|)$ holds for any given $l$ and $S_l$ is measurable with respect to $\sigma \big ( \{ Z_{j+kp_n} | k=0,1,\dots,l \} \big )$. Thus, for any given $\lambda>0$, $S_{l}\in L^{1}\big ( \{ Z_{j+kp_n} | k=0,1,\dots,l \} \big ),\mathbb{P})$. Additionally, 
we obtain from Assumption \ref{as 1}  and $|h|_{SN}\leq B$ that 

\begin{align*}
    |g_j|_{SN}\leq \exp \Big ( {\frac{\lambda A}{|I_j|}} \Big ) \Big |\frac{\lambda}{|I_j|} h \Big |_{SN}\leq \exp \Big ( {\frac{\lambda A}{|I_j|}} \Big ) \frac{\lambda B}{|I_j|}< \infty,    
\end{align*}
for any $ \lambda>0$. Therefore, according to Definition \ref{def 1}, we have, for any $l=1,\dots,|I_j|-1$
\begin{align}
\label{hd4}
    &E[S_{|I_j|-l}]=E\Big [
    S_{|I_j|-l-1}
    \exp\Big (\frac{\lambda h(Z_{j+(|I_j|-l)p_n})}{|I_j|}\Big )\Big ]
    \leq R_1 + R_2 
        \end{align}
        where
    \begin{align*}
   R_1 &=  \left|\text{Cov}\left( S_{|I_j|-l-1} ,\exp\left(\frac{\lambda h(Z_{j+(|I_j|-l)p_n})}{|I_j|}\right)\right)\right| , \\ 
   R_2 & = E\left[S_{|I_j|-l-1} \right ] E\left[\exp\left(\frac{\lambda h(Z_{j+(|I_j|-l)p_n})}{|I_j|}\right)\right].
\end{align*}
By \eqref{eq 2-7} and Definition \ref{def 1} we obtain for the term 
\begin{align*}
    R_1\leq \nu^{-bp_n^{\gamma}}\left\| g_j 
    \right\|
    E\left[S_{|I_j|-l-1} \right ], 
\end{align*}
where $g_j$ is defined in \eqref{hd2}.
Therefore, we obtain from \eqref{hd4} that 
\begin{align}
\label{eq:A6}
E[S_{|I_j|-l}] 
\leq \left(\nu^{-bp_n^{\gamma}}\left\|
  g_j
  \right\|+E\left[\exp\left(\frac{\lambda h(Z_{j+(|I_j|-l)p_n})}{|I_j|}\right)\right]\right)E\left[S_{|I_j|-l-1}\right].
\end{align}
Applying  this argument recursively yields
\begin{align}
\label{eq:A7}
    E_j  = E [ S_{|I_j|-1} ] &  \leq \prod_{l=1}^{|I_j|-1}\left(\nu^{-bp_n^{\gamma}}\left\|
    g_j 
    \right\|+E\left[\exp\left(\frac{\lambda h(Z_{j+lp_n})}{|I_j|}\right)\right]\right)E[S_0]\notag \\
    &=\prod_{l=1}^{|I_j|-1}\left(\nu^{-bp_n^{\gamma}}\left\|g_j
    \right\|+E\left[\exp\left(\frac{\lambda h(Z_{j+lp_n})}{|I_j|}\right)\right]\right)E\left[\exp\left(\frac{\lambda h(Z_j)}{|I_j|}\right)\right]
\end{align} 
(note that  $Z_j$ is the first element of set $I_j$, $j=1,2,\dots,p_n$). 
Recalling the definition of $g_j$ in \eqref{hd2}, the notation   \eqref{eq 2-7} and Assumption \ref{as 1} we obtain  (using $\lambda<\frac{|I_j|}{A}$ in the two inequalities )
\begin{align}
\label{eq:A8}
    \left\| g_j \right\|&= \left\| g_j \right\|_{\infty}+
    \left|g_j \right|_{SN}
    \leq \left\|g_j
    \right\|_{\infty}\left(1+\frac{|h|_{SN}}{A}\right) \leq e\left(\frac{A+B}{A}\right).
\end{align}
On the other hand, since $E[h(Z_j)]=0$, $E[h^2(Z_j)]\leq \sigma^2$, $||h||_{\infty}\leq A$, for any $l$, we have 
(for $0 < \lambda < {|I_j| \over A}$
\begin{align*}
    \exp\left(\frac{\lambda h(Z_{j+lp_n})}{|I_j|}\right)&\leq 1+\frac{\lambda h(Z_{j+lp_n})}{|I_j|}+\frac{1}{2}\left(\frac{\lambda h(Z_{j+lp_n})}{|I_j|}\right)^{2}+\sum_{q=3}^{\infty}\frac{1}{q!}\left(\frac{\lambda h(Z_{j+lp_n})}{|I_j|}\right)^{q}\\
    &\leq 1+\frac{\lambda h(Z_{j+lp_n})}{|I_j|}+\frac{1}{2}\left(\frac{\lambda h(Z_{j+lp_n})}{|I_j|}\right)^{2}\left(1+\sum_{q=3}^{\infty}\left(\frac{\lambda ||h(Z_{j+lp_n})||_{\infty}}{|I_j|}\right)^{q-2}\right)\\
    &\leq  1+\frac{\lambda h(Z_{j+lp_n})}{|I_j|}+\frac{1}{2}\left(\frac{\lambda h(Z_{j+lp_n})}{|I_j|}\right)^{2}\left(\frac{1}{1-\frac{\lambda A}{|I_j|}}\right)\\
    &= 1+\frac{\lambda h(Z_{j+lp_n})}{|I_j|}+\frac{\lambda^{2}h^{2}(Z_{j+lp_n})}{2|I_j|(|I_j|-\lambda A)},
\end{align*}
which implies
\begin{align}
\label{eq:A9}
    E\left[\exp\left(\frac{\lambda h(Z_{j+lp_n})}{|I_j|}\right)\right]\leq 1+\frac{\lambda^{2}\sigma^{2}}{2|I_j|(|I_j|-\lambda A)}\leq \exp\left(\frac{\lambda^{2}\sigma^{2}}{2|I_j|(|I_j|-\lambda A)}\right).
\end{align}
Using the estimates \eqref{eq:A9} and \eqref{eq:A8} in  \eqref{eq:A7} and letting $p_n=\left(\frac{\omega}{b}\log_{\nu}n\right)^{\frac{1}{\gamma}}$, we thus have
\begin{align}
\label{eq:A10}
    E_{j}&\leq \left(\nu^{-bp_n^{\gamma}}e\left(\frac{A+B}{A}\right)+\exp\left(\frac{\lambda^{2}\sigma^{2}}{2|I_j|(|I_j|-\lambda A)}\right)\right)^{|I_j|-1}E[S_0]\\ \notag
    &\leq e\left(\nu^{-bp_n^{\gamma}}e\left(\frac{A+B}{A}\right)+1\right)^{|I_j|}\exp\left(\frac{\lambda^{2}\sigma^{2}}{2(|I_j|-\lambda A)}\right)\\\notag
    &\leq e\exp\left(e\left(\frac{A+B}{A}\right)\nu^{-bp_n^{\gamma}}|I_j|\right)\exp\left(\frac{\lambda^{2}\sigma^{2}}{2(|I_j|-\lambda A)}\right)\\\notag
    &=e\exp\left(e \nu^{\log_{\nu}|I_j|+\log_{\nu}\left(\frac{A+B}{A}\right)-bp_n^{\gamma}}\right)\exp\left(\frac{\lambda^{2}\sigma^{2}}{2(|I_j|-\lambda A)}\right)\\\notag
    &\leq e\exp\left(e \nu^{\log_{\nu}n+\log_{\nu}\left(\frac{A+B}{A}\right)-bp_n^{\gamma}}\right)\exp\left(\frac{\lambda^{2}\sigma^{2}}{2(|I_j|-\lambda A)}\right)\\\notag
    &\leq e\exp\left(3\left(\frac{A+B}{An^{\omega-1}}\right)\right)\exp\left(\frac{\lambda^{2}\sigma^{2}}{2(|I_j|-\lambda A)}\right)\\\notag
    &\leq e^{2}\exp\left(\frac{\lambda^{2}\sigma^{2}}{2(|I_j|-\lambda A)}\right),
\end{align}
where the fourth inequality uses the fact that $\log_{\nu} |I_j|\leq \log_{\nu} n$ and  the last inequality is a consequence of the condition that $n\geq N_{0}:=\min\{n\in\mathbb{N}, \frac{A+B}{An^{\omega-1}}\leq \frac{1}{3}\}$. Note that the upper bound of $E_j$ is independent from $j$, which yields (observing \eqref{eq a-3})
\begin{align}
    \label{eq:A11}
E\left[\exp\left(\frac{\lambda}{n}\sum_{i=1}^{n}h(Z_i)\right)\right]\leq e^{2}\exp\left(\frac{\lambda^{2}\sigma^{2}}{2(|I_j|-\lambda A)}\right) 
\end{align}
for all $0 < \lambda < {A \over |I_j|} $.
Hence, it follows by Markov's inequality that
\begin{align} 
\nonumber \mathbb{P}\left(\left|\frac{1}{n}\sum_{i=1}^{n}h(Z_i)\right|>t\right) & <2e^{2}\exp\left(\frac{\lambda^{2}\sigma^{2}}{2(|I_j|-\lambda A)} -\lambda t\right) \\
\nonumber
&<2e^2\max_{j}\exp\left(-\frac{|I_j|t^{2}}{2(\sigma^{2}+tA)}\right),
\end{align}
where we have used  $\lambda=\frac{|I_j|t}{\sigma^{2}+tA} <\frac{|I_j|}{A} $ for the second inequality.
Together with the fact that
\begin{align*}
    |I_j|\geq \frac{n}{2p_n}=\frac{n}{2(\frac{\omega}{b}\log_{\nu}n)^{\frac{1}{\gamma}}},
\end{align*}
 the proof of Theorem \ref{Th 1} is completed.
 \hfill $\Box$ \\\\
 \noindent
 {\bf Proof of  Theorem \ref{Th 2}.} The proof uses simialr arguments  as given in the proof of Theorem \ref{Th 1}, and  we only highlight the differences. Using  the same notations as in the proof of Theorem \ref{Th 1} the assumption on the decay of the mixing coefficients, we  obtain  for fixed $j$ and $l=1,...,|I_j|$,
\begin{align*}
E[S_{|I_j|-l}]\leq \left(p_{n}^{-\gamma}\left\|\exp\left(\frac{\lambda h(z)}{|I_j|}\right)\right\|+E\left[\exp\left(\frac{\lambda h(Z_{|I_j|-l})}{|I_j|}\right)\right]\right)E[S_{|I_j|-l-1}].    
\end{align*}
Then, using  this argument recursively, it follows that 
\begin{align*}
    E_{j}\leq \prod_{l=2}^{|I_j|}\left(p_{n}^{-\gamma}\left\|\exp\left(\frac{\lambda h(z)}{|I_j|}\right)\right\|+E\left[\exp\left(\frac{\lambda h(Z_{|I_j|-l})}{|I_j|}\right)\right]\right)E[S_1].
\end{align*}
Note that, by assuming $\lambda<\frac{|I_j|}{A}$, 
\begin{align*}
    \left\|\exp\left(\frac{\lambda h(z)}{|I_j|}\right)\right\|&\leq \left\|\exp\left(\frac{\lambda h(z)}{|I_j|}\right)\right\|_{\infty}\left(\left|\frac{\lambda h(z)}{|I_j|}\right|_{SN}+1\right)\\
    &\leq \exp\left(\frac{\lambda A}{|I_j|}\right)\frac{A+B}{A}\\
    &\leq e\frac{A+B}{A},
\end{align*}
and  letting $p_{n}=n^{\frac{\alpha+1}{\gamma+1}}$ yields for the term $E_j$ in \eqref{hd1}
\begin{align*}
    E_j&\leq \left(p_{n}^{-\gamma}e\left(\frac{A+B}{A}\right)+\exp\left(\frac{\lambda^{2}\sigma^{2}}{2|I_j|(|I_j|-\lambda A)}\right)\right)^{|I_j|}||S_1||_{\infty}\\
    &\leq e\left(p_{n}^{-\gamma}e\left(\frac{A+B}{A}\right)+1\right)^{|I_j|}\exp\left(\frac{\lambda^{2}\sigma^{2}}{2(|I_j|-\lambda A)}\right)\\
    &\leq e\exp\left(p_{n}^{-\gamma}e\left(\frac{A+B}{A}\right)|I_j|\right)\exp\left(\frac{\lambda^{2}\sigma^{2}}{2(|I_j|-\lambda A)}\right)\\
    &\leq e\exp\left(p_{n}^{-(\gamma+1)}e\left(\frac{A+B}{A}\right)n\right)\exp\left(\frac{\lambda^{2}\sigma^{2}}{2(|I_j|-\lambda A)}\right)\\
    &\leq e\exp\left(en^{-(\alpha+1)}\left(1+\frac{n^{\alpha}}{A}\right)n\right)\exp\left(\frac{\lambda^{2}\sigma^{2}}{2(|I_j|-\lambda A)}\right)\\
    &\leq e\exp\left(e\left(\frac{1}{n^{\alpha}}+\frac{1}{A}\right)\right)\exp\left(\frac{\lambda^{2}\sigma^{2}}{2(|I_j|-\lambda A)}\right).
\end{align*}
Then, using  the condition $n\geq N_{1}:=\min\{n\in\mathbb{N}:n^{\alpha}>A$\}, we obtain
\begin{align*}
    E_{j}\leq \exp\left(1+\frac{2e}{A}\right)\exp\left(\frac{\lambda^{2}\sigma^{2}}{2(|I_j|-\lambda A)}\right).
\end{align*}
We now use the same arguments as in the proof of Theorem \ref{Th 1} to complete the proof of Theorem \ref{Th 2}.
 \hfill $\Box$ 
\\


\noindent
{\bf Proof of Theorem \ref{thm 3}.} The proof is performed in several steps.
\smallskip

{\par \textbf{Step 1:  (Shadow sample and basic inequality).} 
With the notations  \eqref{eq a-1}, \eqref{eq a-2} and the fact that $\sum_{j=1}^{p_n}w_{j}=1$, we  have 
\begin{align}
    \label{eq B1}
   \mathbb{P}\left(\left|\frac{1}{n}\sum_{i=1}^{n}h(Z_{i})\right|>t\right)&\leq \mathbb{P}\left(\sum_{j=1}^{p_n}w_{j}\left|\frac{1}{|I_j|}\sum_{k=0}^{|I_j|-1}h(Z_{j+kp_n})\right|>\sum_{j=1}^{p_n}w_j t\right) \leq \sum_{j=1}^{p_n}  \mathbb{P}_{j} , 
\end{align}
where 
    \begin{align} \label{hd11}
   \mathbb{P}_{j} &= \sum_{j=1}^{p_{n}}\mathbb{P}\left(\left|\frac{1}{|I_j|}\sum_{k=0}^{|I_j|-1}h(Z_{j+kp_n})\right|>t\right) .  
\end{align}
Hence, it suffice to bound each $\mathbb{P}_{j}$. For each $j$, we introduce a group of shadow samples $\{ Z_{j+kp_n}':k=0,1,\dots,|I_j|-1\} $ 
\footnote{This introduction can always be done by enlarging the original probabilistic space, like adding product-type structure.} satisfying 
\begin{itemize}
    \item [(1)] The distributions of  $Z_{j+kp_n}'$ and $Z_{j+kp_n}$ are identical, $k=0,1,\dots,|I_j|-1$.
    \item [(2)] The random variables  $\{ Z_{j+kp_n}':k=0,1,\dots,|I_j|-1\} $ are mutually independent.
\end{itemize}
 Then, by denoting $B_n=\sqrt{\sum_{k=0}^{|I_j|-1}\text{Var}(h(Z_{j+kp_n}))}$, we obtain
 \begin{align}
     \mathbb{P}_{j}&=\mathbb{P}\left(\left|{1 \over |I_j|}\sum_{k=0}^{|I_j|-1}h(Z_{j+kp_n})\right|>t\right)=\triangle_{j}(t)+\mathbb{P}\left(\left|\frac{1}{|I_j|}\sum_{k=0}^{|I_j|-1}h(Z'_{j+kp_n})\right|>t\right) , 
     \label{hd10}
        \end{align}
        where 
        \begin{align} \label{hd12}
\triangle_{j}(t) &:=\mathbb{P}\left(\left|
\sum_{k=0}^{|I_j|-1}h(Z_{j+kp_n})\right|>
{|I_j|}
t\right)-\mathbb{P}\left(\left|
\sum_{k=0}^{|I_j|-1}h(Z'_{j+kp_n})\right|>
{|I_j|}
t\right).
 \end{align}
Next, we introduce the notations 
\begin{align*}
T_{n} &=\sum_{k=0}^{|I_j|-1}h(Z_{j+kp_n}) ~,~~
T_{n}' = \sum_{k=0}^{|I_j|-1}h(Z'_{j+kp_n}) \\
F_{n}(x) &= \mathbb{P}(T_{n}\leq x)  ~,~~
F'_{n}(x) = \mathbb{P}(T_{n}'\leq x) \\
\phi_{n}(u)&=E[e^{iuT_n}]  =\int_{\mathbb{R}}e^{iux}dF_{n}(x) ~,~~
\phi'_{n}(u) =E[e^{iuT_n'}]  =\int_{\mathbb{R}}e^{iux}dF'_{n}(x).
\end{align*}
and obtain 
\begin{align*}
    |\triangle_{j}(t)|&=\left|\mathbb{P}\left(|T_{n}|\leq |I_j|t\right)-\mathbb{P}\left(|T'_{n}|\leq |I_j|t\right)\right| \leq 2\sup_{x}|F_{n}(x)-F'_{n}(x)|.
\end{align*}
By  Theorem 1 in Chapter 5 of  \cite{petrov2012sums}, it follows that 
\begin{align*}
    |\triangle_{j}(t)|\leq \triangle_{1T}+\triangle_{2T},
\end{align*}
where $T>0$ and 
\begin{align}
\label{hd13}
    \triangle_{1T}&=2b_{0}\int_{-T}^{T}\left|\frac{\phi_{n}(u)-\phi'_{n}(u)}{u}\right|du\\
    \triangle_{2T}&=2b_0 T\sup_{x}\int_{|y|\leq \frac{C(b_0)}{T}}|F_{n}'(x+y)-F_{n}'(x)|dy , 
    \label{hd14}
\end{align}
where the constants $b_{0}$ and $C(b_0)$ are defined in \cite{petrov2012sums}.
\\


\textbf{Step 2: (bound for $\triangle_{2T}$).} Note that
\begin{align*}
 \triangle_{2T}&=2b_0 T\sup_{x}\int_{0}^{\frac{C(b_0)}{T}}(Q'_{n}(x+y)-Q'_{n}(x))dy+2b_0 T\sup_{x}\int_{-\frac{C(b_0)}{T}}^{0}(Q'_{n}(x)-Q'_{n}(x+y))dy\\
    & \leq 2b_0 T\sup_{x}\int_{0}^{\frac{C(b_0)}{T}}(Q'_{n}(x+y)-Q'_{n}(x))dy + 2b_0 T\int_{-\frac{C(b_0)}{T}}^{0}\sup_{x}(Q'_{n}(x)-Q'_{n}(x+y))dy \\
        & = 2b_0 T\sup_{x}\int_{0}^{\frac{C(b_0)}{T}}(Q'_{n}(x+y)-Q'_{n}(x))dy + 2b_0 T\int_{-\frac{C(b_0)}{T}}^{0}\sup_{x}(Q'_{n}(x-y)-Q'_{n}(x))dy \\
    & =  4b_0T\int_{0}^{\frac{C(b_0)}{T}}\sup_{x}(Q'_{n}(x+y)-Q'_{n}(x))dy \\
    & \leq4b_0T\int_{0}^{\frac{C(b_0)}{T}}\mathbf{Q}(W'_{n};y)dy ,
\end{align*}
where the last equality follows by a change of variable  and
$$
\mathbf{Q}(W'_{n};y) = \sup_x \mathbb{P} (x \leq W_n' \leq x +y) 
$$
denotes the {\it concentration function} of the random variable $W_n'$,
and we have used the fact that
\begin{align*}
    2b_0 T\int_{-\frac{C(b_0)}{T}}^{0}\sup_{x}(Q'_{n}(x)-Q'_{n}(x+y))dy
   & =2b_0 T\int_{-\frac{C(b_0)}{T}}^{0}\sup_{x}(Q'_{n}(x-y)-Q'_{n}(x))dy\\
    & = 2b_0T\int_{0}^{\frac{C(b_0)}{T}}\sup_{x}(Q'_{n}(x+y)-Q'_{n}(x))dy.
\end{align*}
Therefore, Lemma 3 on page 38 in \cite{petrov2012sums} with $a=\frac{T}{C(b_0)}$, $T=n^{\delta+1\lor\beta}C(b_0)$ and Assumption \ref{as 3}
yield for any $\delta>0$
\begin{align*}
    \triangle_{2T}&\leq 4b_0T\Big(\frac{96}{95}\Big)^2\Big(\int_{-a}^{a}\prod_{i=1}^n|\psi_{h(Z_i)}(u)|du\Big)\int_{0}^{\frac{C(b_0)}{T}}\max\{y,\frac{1}{a}\}dy\\
    &\leq 4b_0T\Big(\frac{96}{95}\Big)^2\Big(\int_{-a}^{a}\prod_{i\in I'_n}|\psi_{h(Z_i)}(u)|du\Big)\int_{0}^{\frac{C(b_0)}{T}}\max\{y,\frac{1}{a}\}dy\\
     &\leq 4b_0T\Big(\frac{96}{95}\Big)^2\Big(2M_n+\int_{M_n\leq |u|\leq a}(1-\tau_{n})^{c|I_j|}du\Big)\int_{0}^{\frac{C(b_0)}{T}}\max\{y,\frac{1}{a}\}dy\\
    &=\frac{4b_0C^2(b_0)}{T}\Big(\frac{96}{95}\Big)^2\Big(2M_n+2(n^{\delta+1\lor\beta}-M_n)(1-\tau_{n})^{|I_j|}\Big).
\end{align*}
Let $\Theta_\delta$ denote a constant specified below. It depends on $\delta>0$ 
but not on the sample size $n$. Define $p_n=(b^{-1}\Theta_\delta\log n)^{1/\gamma}$. Then, there exists some universal $C>0$ such that 
\begin{align*}
    \triangle_{2T}&\leq C\Big(\frac{M_n}{n^{\delta+1\lor\beta}}+n^{\delta+(1\lor \beta)}(1-\tau_{n})^{|I_j|}\Big) \\
    &\leq C\Big(\frac{M_n}{n^{\delta+(1\lor \beta)}}+n^{\delta+(1\lor \beta)}(1-\tau_{n})^{\frac{n}{2p_n}}\Big)\\
   & \leq 
   C\Big(\frac{M_n}{n^{\delta+(1\lor \beta)}}+n^{\delta+(1\lor \beta)}
 e^{-\frac{n\tau_{n}}{2p_n}}
  \Big) \\
  & =    C\Big(\frac{M_n}{n^{\delta+(1\lor \beta)}}+n^{\delta+(1\lor \beta)}
 e^{-\frac{n\tau_{n}}{2(b^{-1}\Theta_\delta\log n)^{1/\gamma}}}
  \Big) , 
\end{align*}
where we have used the inequality $1-x \leq e^{-x}$.
By the conditions on
 $\tau_n$ and $M_n$ in Assumption \ref{as 3}, there exists a universal $C_{A}>0$ such that  the inequality 
\begin{align*}
    \triangle'_{2T}\leq C_{A}(\frac{1}{n^{\delta+(1\lor \beta)-\beta}}+\frac{1}{n^{\kappa-\delta-(1\lor \beta)}})
\end{align*}
holds for any given $\kappa>0$ and $\delta>0$. Thus, by letting $\kappa=2(\delta+(1\lor \beta))$, we have
\begin{align}
\label{eq B5}
     \triangle'_{2T}\leq C_{A}(\frac{1}{n^{\delta}}+\frac{1}{n^{\delta}})= \frac{C'_A}{n^{\delta}},\ C'_{A}=2C_A.
\end{align}
\par \textbf{Step 3: (bound for $\triangle_{1T}$).} 
 For a real sequence  $(\epsilon_n)_{n \in \mathbb{N} }$ with  $\lim_{n \to \infty }\epsilon_{n}\searrow 0$, we use the decomposition 
\begin{align*}
    \triangle_{1T}
    &=:\triangle^1_{1T}+\triangle^2_{1T}
\end{align*}
where 
    \begin{align}
    \nonumber 
   \triangle^1_{1T}=&2b_0\int_{\epsilon_{n}\leq |u|\leq T}\left|\frac{\phi_{n}(u)-\phi'_{n}(u)}{u}\right|du, \\
   \triangle^2_{1T}=&
   2b_0\int_{-\epsilon_{n}}^{\epsilon_{n}}\left|\frac{\phi_{n}(u)-\phi'_{n}(u)}{u}\right|du~.
   \label{hd15}
\end{align}
We will derive bounds for $\triangle_{1T}^1$ and $\triangle_{1T}^2$ separately by delivering different upper bounds to $|\phi_{n}(u)-\phi_{n}'(u)|$.
 To bound $\triangle^{1}_{1T}$, recall that, for each $j$, the distribution of $Z_{j+kp_n}$ and $Z_{j+kp_n}'$ are identical but the random variables $Z_{j+kp_n}'$, $k=0,1,\dots,|I_j|-1$, are mutually independent. Thus, by using the notations $W_{j+kp_n}= h(Z_{j+kp_n})$ and $W_{j+kp_n}'=h(Z_{j+kp_n}')$, 
we have
\begin{align*}
    \Psi(u):=&|\phi_{n}(u)-\phi_{n}'(u)|=|E[e^{iuT_{n}}]-E[e^{iuT_n'}]|\\
    =& \Big |E \Big  [e^{iu\sum_{k=0}^{|I_j|-1}W_{j+kp_n}} \Big  ]-E[e^{iu\sum_{k=0}^{|I_j|-1}W_{j+kp_n}'}] \Big |\\
    =&\left|E[\prod_{k=0}^{|I_j|-1}e^{iuW_{j+kp_n}}]-\prod_{k=0}^{|I_j|-1}E[e^{iuW_{j+kp_n}}]\right|  .
\end{align*}
This gives the estimate 
\begin{align*}
    \Psi(u)& \leq  \sum_{l=1}^{|I_j|-1}  \Big |E \Big  [\prod_{k=0}^{|I_j|-l-1}e^{iuW_{j+kp_n}}e^{iuW_{j+(|I_j|-l)p_n}}\Big ]-E \Big [\prod_{k=0}^{|I_j|-l-1}e^{iuW_{j+kp_n}}\Big  ]E \Big [e^{iuW_{j+( |I_j|-l)p_n}}\Big  ]\Big  |\Xi_{l},
\end{align*}
where
\begin{align*}
    \Xi_{l}&=\left\{
\begin{aligned}
&1, & \text{when $l=1$}, \\
&\prod_{k=|I_j|-l+1}^{|I_j|-1}\left|E[e^{iuW_{j+kp_n}}]\right|, & \text{when $2\leq l\leq {|I_j|-1}$}.
\end{aligned}\right.
\end{align*}
As  $0\leq \Xi_{l}\leq 1$  ($l=1,2,\dots,|I_j|-1$), we obtain

\begin{align}
    \Psi(u)\leq  \sum_{l=1}^{|I_j|-1}E_{l},
    \label{hd8}
\end{align}
where
$$
E_l:=
\sum_{l=1}^{|I_j|-1} \Big |E \Big [\prod_{k=0}^{|I_j|-l-1}e^{iuW_{j+kp_n}}e^{iuW_{j+( |I_j|-l)p_n}}\Big ]-E\Big [\prod_{k=0}^{|I_j|-l-1}e^{iuW_{j+kp_n}} \Big ]E \Big [e^{iuW_{j+( |I_j|-l)p_n}}\Big ] \Big |.
$$
We estimate  each $E_l$ by 
\begin{align}
  \label{hd9}
    E_{l}=&\left|\text{Cov}\left(e^{iuT_{|I_j|-l-1}},e^{iuW_{j+( |I_j|-l)p_n}}\right)\right| \leq E_{l1}+E_{l2}+E_{l3}+E_{l4} ,
    \end{align} 
    where 
    \begin{align*}
 E_{l1}  = & \big |\text{Cov}\big ( \cos (uT_{|I_j|-l-1}),\cos (W_{j+( |I_j|-l)p_n} ) \big  )|,~ E_{l2} =  \big  |\text{Cov} \big  (\cos (uT_{|I_j|-l-1}),\sin ( W_{j+( |I_j|-l)p_n}) \big  )\big  |\\
    E_{l3}= &\big  |\text{Cov}\big  (\sin (uT_{|I_j|-l-1}),\cos ( W_{j+( |I_j|-l)p_n}) \big  )\big  |,~E_{l4} = \big |\text{Cov}\big  (\sin (uT_{|I_j|-l-1}),\sin (W_{j+( |I_j|-l)p_n}) \big )\big  |
\end{align*}

For $E_{l1}$, according to Assumption \ref{as 2} and the assumption that $|h|_{SN}\leq n^{B}$, we have
\begin{align*}
    \left|\cos (uh)\right|_{SN}\leq C_{B}|u| ^{\eta_0}|h|_{SN}\leq C_{B}|u|^{\eta_0}n^{B},
\end{align*}
which implies 
\begin{align*}
    E_{1l}&\leq C(p_n)||\cos (uT_{|I_j|-l-1})||_{\infty}||\cos (uh)||\\
    &\leq C(p_n)||\cos (uT_{|I_j|-l-1})||_{\infty}(||\cos (uh)||_{\infty}+ \left|\cos (uh)\right|_{SN})\\
    &\leq C(p_n)(1+ C_{B}|u|^{\eta_0}n^{B}).
\end{align*}
By the same argument we get the same bounds for $E_{l2}, E_{l3}$ and $E_{l4}$  and we 
obtain from \eqref{hd8} and \eqref{hd9} that 
the inequality  
\begin{align}
\label{eq B6}
    \Psi(u) \leq 
    4|I_j|C(p_n)(1+ C_{B}|u|^{\eta_0}n^{B})
\end{align}
holds for every $n\in\mathbb{N}$ and $u\in\mathbb{R}$. Then, together with the definition of $\triangle_{1T}^1$ it follows that 
\begin{align}
\label{eq triangle_{1T}^1}
    \triangle_{1T}^{1}\leq 16b_0C(p_n)\left(n\left(\log T+\log \frac{1}{\epsilon_n}\right)+C_B\eta_0^{-1}n^{B+1}T^{\eta_0}\right),
\end{align}
for any $T>0$ and $\epsilon_n\searrow 0$. 

We now continue deriving an upper bound for  $\triangle^2_{1T}$, noting  that
\begin{align*}
    \Psi(u)=|\phi_n(u)-\phi_{n}'(u)|\leq |E[\cos uT_{n}]-E[\cos uT'_{n}]|+|E[\sin uT_{n}]-E[\sin uT'_{n}]|.
\end{align*}
With the notation 
\begin{align*}
    S_{k}=W_{j+p_n}+\dots+W_{j+(k-1)p_n} +W'_{j+(k+1)p_n}+\dots+W'_{j+(|I_j|-1)p_n},
\end{align*}
we  obtain the representations
\begin{align*}
    T_n  & =S_{|I_j|-1}+W_{j+(|I_j|-1)p_n},
    \\ 
    T'_{n}& =S_0+W'_{j+p_n}, \\ 
    S_{k-1}+W_{j+(k-1)p_n}  & =S_{k}+W'_{j+kp_n} ~~~~ ~~ (k=1, \ldots   , |I_j|-2).
\end{align*}
If $f:\mathbb{R}\rightarrow\mathbb{R}$ is a continuously differentiable function 
such that $\| f'\|_\infty < \infty $, we have
\begin{align*}
    |E[f(uT_n)]-E[f(uT'_n)]|& \leq \sum_{k=0}^{|I_j|-1}|E[f(u(S_k+W_{j+kp_n}))]-E[f(u(S_k+W'_{jk}))]| \\
    & \leq \sum_{k=0}^{|I_j|-1}  2||f'||_{\infty}|u|E|W_{j+kp_n}| \\
    &=  \sum_{k=0}^{|I_j|-1}2||f'||_{\infty}|u|E|f(Z_{j+kp_n})|,  
\end{align*}
where we have used the mean-value theorem, that is 
\begin{align*}
    f(u(S_k+W_{j+kp_n}))&=f(uS_k)+f'(u(S_{k}+\xi_{jk}))uW_{j+kp_n},\\
     f(u(S_k+W'_{jk}))&=f(uS_k)+f'(u(S_{k}+\xi'_{jk}))uW'_{jk},
\end{align*}
for some random variables  $\xi_{jk}$ and  $\xi'_{jk}$.
This yields
\begin{align}
    \label{eq B7}
    \Psi(u) =|\phi_n(u)-\phi_{n}'(u)|  \leq 4|u|\sum_{k=0}^{|I_j|-1}E|h(Z_{j+kp_n})|\leq 4|I_j||u| \max_{i}E|h(Z_i)|.
\end{align}
Observing the  definition of $\triangle_{1T}^2$ in \eqref{hd15} it follows that
\begin{align}
    \label{triangle_1T^2}
    \triangle_{1T}^2&\leq 4|I_j|(\max_{i}E|h(Z_i)|)\int_{-\epsilon_n}^{\epsilon_n}du\leq 8n(\max_{i}E|h(Z_i)|)\epsilon_n
\end{align}
Therefore, 
\begin{align}
    \label{eq B8}
    \triangle_{1T}\leq 8n(\max_{i}E|h(Z_i)|)\epsilon_n+ 16b_0C(p_n)\left(n\left(\log T+\log \frac{1}{\epsilon_n}\right)+C_B\eta_0^{-1}n^{B+1}T^{\eta_0}\right),
\end{align}
and combining  \eqref{eq B5} and \eqref{eq B8} yields 
\begin{align}
    \label{eq B9}
    \max_{1\leq j\leq p_n}\sup_{t>0}|\triangle_{j}(t)|\leq& \frac{C'_A}{n^{\delta}}+8n(\max_{i}E|h(Z_i)|)\epsilon_n\notag\\
    &+ 16b_0C(p_n)\left(n\left(\log T+\log \frac{1}{\epsilon_n}\right)+C_B\eta_0^{-1}n^{B+1}T^{\eta_0}\right)
\end{align}
 for any $T>0$, $\epsilon_{n}\searrow 0$ and $p_n>0$. Recall we have let $T=C(b_0)n^{\delta+(1\lor \beta)}$ and and $p_n=(b^{-1}\Theta_{\delta}\log_{\nu}n)^{\frac{1}{\gamma}}$. By setting $\epsilon_{n}=\frac{1}{n^{\delta+2}}$, we obtain from \eqref{eq B9} that 
\begin{align*}
    \max_{1\leq j\leq p_n}\sup_{t>0}|\triangle_{j}(t)| \leq &\frac{C'_A}{n^{\delta}}+\frac{8(\max_{i}E|h(Z_i)|)}{n^{\delta+1}}\\&+\frac{32b_0(\delta+\beta\lor 2)n\log( (C(b_0)\lor 1) n)}{n^{\Theta_{\delta}}}\\&+\frac{16b_0C_B(C(b_0))^{\eta_0}\eta_0^{-1}n^{B+1+\eta_0(\delta+1\lor\beta)}}{n^{\Theta_{\delta}}}.
\end{align*}
By setting
\begin{align}
\label{1}
    \Theta_{\delta}=B+2+(\eta_0+1)(\delta+1\lor\beta),
\end{align}
we obtain the inequality 
\begin{align}
    \label{eq B10}
    &\max_{1\leq j\leq p_n}\sup_{t>0}|\triangle_{j}(t)|\leq \frac{C_{\delta}}{n^{\delta}},
\end{align}
with the constants 
\begin{align}    
    \label{2}
 C_{\delta}&= C'_A+8\max E|h(Z_i)|+32b_0(\delta+\beta\lor 2)+16b_0C_B(C(b_0))^{\eta_0}\eta_0^{-1}\notag\\
 &=:C'_{\delta}+8\max E|h(Z_i)|.
\end{align}

  \textbf{Step 4}
Combining \eqref{eq B1}, \eqref{hd10} and  \eqref{eq B10}  with Bernstein inequality for independent random variables, we obtain
for the right hand side of \eqref{eq B1} that
\begin{align*}
    \sum_{j=1}^{p_n}\mathbb{P}_j&\leq 2p_n\exp\left(-\frac{nt^2}{p_n(\sigma^{2}_{\max}+\frac{At}{3})}\right)+\frac{p_nC_{\delta}}{n^{\delta}} \\ 
    &\leq2p_n\exp\left(-\frac{nt^2}{p_n(\sigma^{2}_{\max}+\frac{At}{3})}\right)+\frac{C’_{\delta}+8\max_{i} E|h(Z_i)|}{n^{\delta-1}}, 
\end{align*}
where $\sigma^2_{\max}=\max_{1\leq j\leq p_n}\frac{1}{|I_j|}\sum_{k=0}^{|I_j|-1}\text{Var}(h(Z_{j+kp_n}))$
and the last inequality follows from \eqref{2} by choosing a constant $\delta>1$.
} 
\hfill $\Box$
\\ \\
{\bf Proof of Proposition \ref{prop A2}}{
We denote $(X_1,Y_1)$ and $(x,y)$ as $Z=(Z_1,...,Z_D)$ and $z=(z_1,...,z_D)$ respectively and use the representation 
 $$
 K_{x,h}(X_1)K_{y,h}(Y_1) =\prod_{s=1}^{D}K\Big (\frac{Z_s-z_s}{h} \Big )=:g(Z).
 $$
Since $K$ has compact support $[-1,1]$ (Assumption \ref{kernel}) we obtain 
\begin{align*}
   \psi_{g(Z)}(u)=\int_{B_{D}(z,h)}e^{iu\prod_{s=1}^{D}K(\frac{w_s-z_s}{h})}f_{(X,Y)}(w)dw+1-\mathbb{P}^{(X,Y)}
    (B_{D}(z,h)) ,
\end{align*}
where  $B_{D}(z,h):=\prod_{s=1}^{D}[z_s-h,z_s+h]$.
By the  smoothness assumptions on  $f_{(X,Y)}$ we have 
\begin{align*}
  \mathbb{P}^{(X,Y)}(B_{D}(z,h))&=\int_{B_{D}(z,h)}|f_{(X,Y)}(w)-f_{(X,Y)}(z)+f_{(X,Y)}(z)|dw\\
  &\geq \int_{B_{D}(z,h)}f_{(X,Y)}(z)dw-\int_{B_{D}(z,h)}|f_{(X,Y)}(w)-f_{(X,Y)}(z)|dw\\
  &\geq (2h)^{D}f_{(X,Y)}(z)-D2^{D}C_{f_{(X,Y)}}h^{\min\{\alpha,1\}+D},
\end{align*}
where $C_{f_{(X,Y)}}$ is the Hölder constant of the  function $f_{(X,Y)}$. 
Consequently, there exist constants $C_1,C_1'>0 $ (independent of the bandwidth $h$ and sample size $n$) such that
\begin{align}
\label{det5}
    |\psi_{h(Z)}(u)|&\leq 1-C_1h^D+ C_1'h^{\min\{\alpha,1\}+D}+E(u),
\end{align}
where 
\begin{align*}
   E(u)&= \Big |\int_{B_{D}(z,h)}e^{iu\prod_{s=1}^{D}K(\frac{w_s-z_s}{h})}f_{(X,Y)}(w)dw\Big |  \leq E_{1}(u)+E_{2}(u) 
\end{align*}
and
\begin{align*}
    E_1(u)&:=  \Big |\int_{B_{D}(z,h)}e^{iu\prod_{s=1}^{D}K(\frac{w_s-z_s}{h})}(f_{(X,Y)}(w)-f_{(X,Y)}(z))dw\Big | , \\
    E_2(u)&:= f_{(X,Y)}(z)\left|\int_{B_{D}(z,h)}e^{iu\prod_{s=1}^{D}K(\frac{w_s-z_s}{h})}dw\right|.
\end{align*}
A simple calculation shows  that there exists some $C_1''>0$ independent of $h$ and $n$ such that
\begin{align}
    E_{1}(u)
    &\leq h^{D}\int_{[-1,1]^D}|f_{(X,Y)}(z+vh)-f_{(X,Y)}(z)|dv\leq C_1''h^{\alpha+D}.
    \label{det5a}
\end{align}
To derive an inequality for  $E_2(u)$ we use Assumption \ref{kernel} and a change of variable and obtain
\begin{align} \nonumber 
    E_2(u)&\leq ||f_{(X,Y)}||_{\infty}(2h)^D\left|\int_{[0,1]^D}e^{iu\prod_{s=1}^{D}K(w_s)}dw_1\dots dw_D\right|\\
    &\lesssim (\epsilon_nh)^D +h^D\int_{[\epsilon_n,1-\epsilon_n]^{D-1}}\left|\int_{\epsilon_{n}}^{1-\epsilon_n}e^{iuK(w_1)\prod_{s=2}^{D}K(w_s)}dw_1\right|dw_2\ldots dw_D.
    \label{det5b} 
\end{align}
for any given $\epsilon_n\searrow 0$. 
  
We now focus on the inner integral and not that by 
 Assumption \ref{kernel} the kernel  $K$ is strictly decreasing on $[0,1]$ and twice differentiable.
 Moreover, by condition (K1) in   
 Assumption \ref{kernel},
 we have 
 %
$\prod_{s=2}^{D}K(w_s)\geq K(1-\epsilon_n)^{D-1}\geq  \epsilon_{n}^{\beta_1(D-1)}$ for 
 $w_2, \ldots , w_D\in [\epsilon_n,1-\epsilon_n]$. Then, by letting $v_1=K(w_1)$, it follows that 
\begin{align}
\nonumber 
    \Big |\int_{\epsilon_{n}}^{1-\epsilon_n}e^{iuK(w_1)\prod_{s=2}^{D}K(w_s)}dw_1 \Big | & \leq\Big |\int_{K(\epsilon_n)}^{K(1-\epsilon_n)}\exp \Big  (iuv_1\prod_{s=2}^{D}K(w_s) \Big  )(K'(K^{-1}(v_1)))^{-1}dv_1\Big |\\
  &   \lesssim |u|^{-1}(\epsilon_{n})^{-[(D-1)\beta_1]}|E_2'(u)|.
  \label{det5c}
\end{align}
where
\begin{align*}
  E_2'(u)   & := \int_{K(\epsilon_n)}^{K(1-\epsilon_n)}\exp \Big  (iuv_1\prod_{s=2}^{D}K(w_s) \Big  )(K'(K^{-1}(v_1)))^{-1}d(iu\prod_{s=2}^{D}K(w_s)v_1)  \\
  &= \frac{\exp(iuv_1\prod_{s=2}^{D}K(w_s))}{K'(K^{-1}(v_1))}\Big|_{K(\epsilon_n)}^{K(1-\epsilon_n)}-\int_{K(\epsilon_n)}^{K(1-\epsilon_n)}\frac{K{''}(K^{-1}(v_1))\exp(iuv_1\prod_{s=2}^{D}K(w_s))}{K'(K^{-1}(v_1))^{3}}dv_1.
\end{align*}
%
and the last identity follows by  integration by part.
Assumption \ref{kernel} yields
\begin{align*}
    |E_2'(u)|\lesssim \frac{1}{|K'(\epsilon_n)|}+\frac{||K''||_{\infty}\overline{K}}{|K'(\epsilon_n)|^{3}}\lesssim\epsilon_{n}^{-3{\beta_2}}.
\end{align*}

As $(\epsilon_n)_{n \in \mathbb{N}}$ can be any vanishing sequence converging to $0$, we can choose its convergence rate  such  that there exists a sequence  $(M_n)_{n \in \mathbb{N}}$ with $M_n\lesssim n$  and  $M_n\epsilon_n^{-3\beta_2}\nearrow \infty$. Then, for every $u\geq M_n$, we have $|E_2(u)|=o(h^{D})$. 
Observing \eqref{det5a},  it follows that  $|E(u)|=o(h^{D})$ and 
the inequality \eqref{det5} shows that, by the condition  on $h$,  that $g(Z)$ satisfies Assumption \ref{as 3}.
\bigskip

{\bf Proof of Theorem \ref{Th 4}.} The proof is divided into two steps. In the first step, we show the uniform convergence rates of the joint density estimator \eqref{det19}. Based on this rate, we show the convergence rate of the Hausdorff loss of our mode set estimator \eqref{eq 4-8} in the second step.}
\\

\textbf{Step 1}
 According to Proposition \ref{prop A2}, Assumptions \ref{as 5}, \ref{kernel} and \ref{kernel 2} guarantee that Theorem \ref{thm 3} can be applied to investigate the uniform convergence rate of joint density estimator. 
 More specifically, using Theorem \ref{thm 3} with $h=\big (\frac{\log n}{n}\big  )^{{1}/({2\alpha+D})}$, we obtain, for any  interior point $(x,y)\in\mathbb{X}\times\mathbb{Y}$ that there exist constants $M^*$, $C_M$, $M_2$ and $M_{\delta}$
 independent of $n$ and $(x,y)$  
such that the  inequality 
\begin{align*}   &\mathbb{P}\left(\left|\frac{1}{nh^{D}}\sum_{i=1}^{n}K_{x,h}(X_i)K_{y,h}(Y_i)-f_{(X,Y)}(x,y)\right|>M\left(\frac{\log n}{n}\right)^{\frac{\alpha}{(2\alpha+D)}}\right)\leq \frac{M_2}{n^{C_M}}+\frac{M_{\delta}}{n^{\delta-1}}.
 \end{align*}
 holds for any $M>M^*$, $\delta>1$ 
 and $n\geq N_0$
  (this constant is defined in Theorem \ref{thm 3}).
  Moreover, the constant $C_M$ is increasing with  $M$.
Then, using the same arguments as the proof of Theorem 2 in \cite{hansen2008uniform},
it follows that 
    \begin{align}
    \label{eq 64}
\mathbb{P}\Big (\limsup_{n \to \infty }\sup_{y\in \mathbb{Y}} \big |\hat{f}(x_0,y)-f(x_0,y) \big |\geq \xi_n\Big (\frac{C\log n}{n} \Big)^{\frac{\alpha}{2\alpha+D}}\Big )=o(1)
    \end{align}
where  $C>1$ is constant  and   $\xi_n $ an arbitrary monotone sequence converging to $\infty $.   Furthermore, 
$h=\big (\frac{\log n}{n}\big  )^{{1}/({2\alpha+D})}$, we can further show, that for any  $\epsilon>0$, there exists sufficiently large $M>0$ and $N\in\mathbb{N}^+$ such that
    \begin{align*}   \mathbb{P}\Big  (\sup_{y\in \mathbb{Y}}|\hat{f}(x_0,y)-f(x_0,y)|>M \Big (\frac{\log n}{n} \Big )^{\frac{\alpha}{2\alpha+D}}\Big )<\epsilon
    \end{align*}
    holds for every $n>N$, which means that 
\begin{align*}
    \sup_{y\in \mathbb{Y}}\big |\hat{f}(x_0,y)-f(x_0,y)\big |=\mathcal{O}_{\mathbb{P}}\Big (\frac{\log n}{n}\Big )^{\frac{\alpha}{2\alpha+D}}.
\end{align*}

\textbf{Step 2} Define the event $A_{n}=\{\text{Card}(\widehat{\mathcal{Y}}_{x_0})=K_{x_0}\}$, then 
\begin{align*}
    A^c_{n}=\{\text{Card}(\widetilde{\mathcal{Y}}_{L,x_0})\neq K_{x_0}\}\subset S_{n} := \Big \{\sup_{y\in\mathbb{Y}}|\hat{f}(x_0,y)-f(x_0,y)|>0 \Big \},
\end{align*}
From \eqref{eq 64} it follows that $\lim_{n \to \infty}\mathbb{P}(S_n)=0$ which implies 
  $\lim_{n \to \infty }P(A_n^c)=0$. 
 %
Consequently, we obtain
\begin{align*}
    \mathbb{P}\Big(\sup_{y\in\widehat{\mathcal{Y}}_{x_0}}d(y,\mathcal{Y}_{x_0})>t\Big)&\leq P\Big(\{\sup_{y\in\widehat{\mathcal{Y}}_{x_0}}d(y,\mathcal{Y}_{x_0})\}\cap A_n\Big)+P(A_n^c).
\end{align*}
Note that, for every $\omega\in A_n$, $\widehat{\mathcal{Y}}_{x_0}(\omega)$ can be written as $\{\hat{y}_{x_0,k}:k=1,2,\dots,K_{x_0}\}$, where $\hat{y}_{x_0,k}$ indicates the $k$-th local maximizer of $\hat{f}(x_0,y)$. 
This gives  
\begin{align*}
    \mathbb{P}\Big(\sup_{y\in\widehat{\mathcal{Y}}_{x_0}}d(y,\mathcal{Y}_{x_0})\}\cap A_n\Big) &  = \mathbb{P}
    \Big ( \Big \{\max_{1\leq k\leq K_{x_0}}d(\hat{y}_{x_0,k},\mathcal{Y}_{x_0})>t\Big\}\cap A_n \Big) 
    \\
    &  \leq  \mathbb{P} \Big ( \max_{1\leq k\leq K_{x_0}}d(\hat{y}_{x_0,k},\mathcal{Y}_{x_0})>t \Big)  \\
    &
     \leq K \max_{1\leq k\leq K_{x_0}}\mathbb{P}(d(\hat{y}_{x_0,k},\mathcal{Y}_{x_0})>t).
\end{align*}
Let $B_{k}(r)=B(y_{x_0,k},r)$ denote the closed ball with center $y_{x_0,k}$ and radius $r$. According to conditions (C1) and (C2) of Assumption \ref{as 5}, there exists a constant  $\eta>0$ such that $\{B_{k}(\eta):k=1,..,K_{x_0}\}$ is a group of disjoint closed balls and $f(x_0,y_{x_0,k})=\sup_{y\in B_k}f(x_0,y)$. For $k_0 =1, \ldots , K_{x_0} $, define the event 
\begin{align*}
     R_{k_0}=\big  \{ 
    \text{there is a unique}\ l(k_0)\in\{1,2, \ldots ,K_{x_0}\}\ \text{such that}\ \hat{y}_{x_0,k_0}\in B_{l(k_0)}(\delta) \big \}, 
\end{align*}
and note that $R^c_{k_0}\subset S_n$. Similar to previous arguments  we have
\begin{align*}
    \mathbb{P}(d(\hat{y}_{x_0,k_0},\mathcal{Y}_{x_0})>t)&\leq \mathbb{P}(\{d(\hat{y}_{x_0,k_0},\mathcal{Y}_{x_0})>t\}\cap R_{l(k_0)})+\mathbb{P}(R^c_{l(k_0)})\\
    &\leq \mathbb{P}(t<||\hat{y}_{x_0,k_0}-y_{x_0,l(k_0)}||_E \leq \eta)+\mathbb{P}(R^c_{l(k_0)}).
\end{align*}
Together with condition C2 of Assumption \ref{as 5}, we have
\begin{align*}
    \mathbb{P}(d(\hat{y}_{x_0,k_0},\mathcal{Y}_{x_0})>t)  & \leq \mathbb{P}(||\hat{y}_{x_0,k_0)}-y_{x_0,l(k_0)}||_E>t)+\mathbb{P}(R^c_{l(k_0)})\\
    \lesssim& \mathbb{P}(|f(x_0,\hat{y}_{x_0,k_0})-f(x_0,y_{x_0,l(k_0)})|>t^{\beta})+\mathbb{P}(R^c_{l(k_0)})\\
    \leq& 
    \mathbb{P}(|f(x_0,\hat{y}_{x_0,k_0})-\hat{f}(x_0,\hat{y}_{x_0,k_0})|>t^\beta/2) \\
    & ~~~~~~~~~~~~~~~+\mathbb{P}(|\hat{f}(x_0,\hat{y}_{x_0,k_0})-f(x_0,y_{x_0,l(k_0)})|>t^{\beta}/2)+\mathbb{P}(R^c_{l(k_0)}) \\
    \leq& \mathbb{P} \Big (\sup_{y\in B_{k_0}}|f(x_0,y)-\hat{f}(x_0,y)|>t^\beta/2 \Big ) \\
    & ~~~~~~~~~~~~~~~+\mathbb{P}
    \Big (|\sup_{y\in B_{k_0}}\hat{f}(x_0,y)-\sup_{y\in B_{k_0}}f(x_0,y)|>t^{\beta}/2\Big )+\mathbb{P}(R^c_{l(k_0)})\\
    \leq &2\mathbb{P} \Big (\sup_{y\in B_{k_0}}|f(x_0,y)-\hat{f}(x_0,y)|>t^\beta/2 \Big )+\mathbb{P}(R^c_{l(k_0)}).
\end{align*}
Summarizing
, we obtain 
\begin{align*}
    \mathbb{P}\Big(\sup_{y\in\widehat{\mathcal{Y}}_{x_0}}d(y,\mathcal{Y}_{x_0})>t\Big)&\leq \mathbb{P}\Big(\{\sup_{y\in\widehat{\mathcal{Y}}_{x_0}}d(y,\mathcal{Y}_{x_0})>t\}\cap A_n\Big)+\mathbb{P}(A^c_n)\\
    &\leq 2K\max_{1\leq k_0\leq K_{x_0}}\mathbb{P}\Big(\sup_{y\in B_{k_0}}|f(x_0,y)-\hat{f}(x_0,y)|>t^\beta/2 \Big)+K\mathbb{P}(R^c_{l(k_0)})+\mathbb{P}(A^c_n)
\end{align*}
Based on Step 1, we have 
\begin{align}
\label{eq A3}
\sup_{y\in\widehat{\mathcal{Y}}_{x_0}}d(y,\mathcal{Y}_{x_0})=\mathcal{O}_{\mathbb{P}}\left(\left(\frac{\log n}{n}\right)^{\frac{\alpha}{\beta(2\alpha+D)}}\right).
\end{align}
As for the second part in Hausdorff-loss \eqref{eq 4-9}, repeating our arguments above yields  
\begin{align}
\label{eq A4}
\sup_{y\in\mathcal{Y}_{x_0}}d(y,\widehat{\mathcal{Y}}_{x_0})=\mathcal{O}_{\mathbb{P}}\Big (\Big (\frac{\log n}{n}\Big )^{\frac{\alpha}{\beta(2\alpha+D)}}\Big).
\end{align}
By combining \eqref{eq A3} and \eqref{eq A4}, we finish the proof of Theorem \ref{Th 4}.

\hfill $\Box$


\noindent

\vskip 0.2in
\bibliography{sample}

\end{document}